\documentstyle[12pt]{article}

\font\twlgot =eufm10 scaled \magstep1
\font\egtgot =eufm8
\font\sevgot =eufm7
\font\twlmsb =msbm10 scaled \magstep1
\font\egtmsb =msbm8
\font\sevmsb =msbm7

\newfam\gotfam
\def\pgot{\fam\gotfam\twlgot}
\textfont\gotfam\twlgot
\scriptfont\gotfam\egtgot
\scriptscriptfont\gotfam\sevgot
\def\got{\protect\pgot}
\newfam\msbfam
\textfont\msbfam\twlmsb
\scriptfont\msbfam\egtmsb
\scriptscriptfont\msbfam\sevmsb
\def\Bbb{\protect\pBbb}
\def\pBbb{\relax\ifmmode\expandafter\Bb\else\typeout{You cann't use
Bbb in text mode}\fi}
\def\Bb #1{{\fam\msbfam\relax#1}}

\newcommand{\gO}{{\got O}}
\newcommand{\gQ}{{\got T}}
\newcommand{\gU}{{\got U}}
\newcommand{\gE}{{\got E}}
\newcommand{\gA}{{\got A}}

\newcommand{\gd}{{\got d}}

\newcommand{\gS}{{\got S}}

\def\thebibliography#1{\section*{References}\list
  {[\arabic{enumi}]}{\settowidth\labelwidth{#1}\leftmargin\labelwidth
    \advance\leftmargin\labelsep
    \usecounter{enumi}}
    \def\newblock{\hskip .11em plus .33em minus .07em}
    \sloppy\clubpenalty4000\widowpenalty4000
    \sfcode`\.=1000\relax}

\def\op#1{\mathop{\fam0 #1}\limits}

\newcommand{\Ker}{{\rm Ker\,}}
\newcommand{\im}{{\rm Im\,}}
\newcommand{\nm}[1]{|{#1}|}
\newcommand{\beq}{\begin{equation}}
\newcommand{\eeq}{\end{equation}}
\newcommand{\ben}{\begin{eqnarray}}
\newcommand{\een}{\end{eqnarray}}
\newcommand{\be}{\begin{eqnarray*}}
\newcommand{\ee}{\end{eqnarray*}}
\newcommand{\bea}{\begin{eqalph}}
\newcommand{\eea}{\end{eqalph}}
\newcommand{\cA}{{\cal A}}

\newcommand{\cP}{{\cal P}}

\newcommand{\cL}{{\cal L}}
\newcommand{\cJ}{{\cal J}}

\newcommand{\cV}{{\cal V}}
\newcommand{\cE}{{\cal E}}

\newcommand{\cQ}{{\cal Q}}

\newcommand{\cS}{{\cal S}}
\newcommand{\cC}{{\cal C}}
\newcommand{\cO}{{\cal O}}

\newcommand{\bL}{{\bf L}}

\newcommand{\bs}{{\bf s}}

\newcommand{\al}{\alpha}
\newcommand{\vr}{\varrho}

\newcommand{\dl}{\delta}
\newcommand{\la}{\lambda}
\newcommand{\La}{\Lambda}
\newcommand{\f}{\phi}
\newcommand{\om}{\omega}

\newcommand{\m}{\mu}

\newcommand{\g}{\gamma}
\newcommand{\G}{\Gamma}
\newcommand{\th}{\theta}

\newcommand{\vt}{\vartheta}
\newcommand{\vf}{\varphi}
\newcommand{\up}{\upsilon}
\newcommand{\lng}{\langle}
\newcommand{\rng}{\rangle}
\newcommand{\di}{{\rm dim\,}}

\newcommand{\si}{\sigma}
\newcommand{\Si}{\Sigma}
\newcommand{\w}{\wedge}
\newcommand{\wt}{\widetilde}

\newcommand{\ol}{\overline}
\newcommand{\dr}{\partial}
\newcommand{\ar}{\op\longrightarrow}
\newcommand{\ot}{\otimes}
\newcommand{\ap}{\approx}

\let\ssection=\section
\renewcommand{\section}{\setcounter{equation}{0}\ssection}

\newcounter{eqalph}[section]
\newcounter{equationa}[section]
\newcounter{example}[section]
\newcounter{remark}[section]
\newcounter{theorem}[section]
\newcounter{proposition}[section]
\newcounter{lemma}[section]
\newcounter{corollary}[section]
\newcounter{definition}[section]

\def\theremark{\arabic{section}.\arabic{remark}}

\def\thedefinition{\arabic{section}.\arabic{definition}}

\newenvironment{proof}{\noindent 
{\it Proof.}}{$\Box$ \medskip}
\newenvironment{rem}{\refstepcounter{remark}\medskip\noindent{\it
Remark \theremark.}}{\medskip}
\newenvironment{theo}{\refstepcounter{definition} 
\bigskip\noindent{\bf Theorem \thedefinition.} \it}{\medskip}
\newenvironment{prop}{\refstepcounter{definition} 
\bigskip\noindent{\bf Proposition \thedefinition.}\it}{\medskip}
\newenvironment{lem}{\refstepcounter{definition} 
\bigskip\noindent{\bf Lemma \thedefinition.}\it}{\medskip}
\newenvironment{cor}{\refstepcounter{definition} 
\bigskip\noindent{\bf Corollary \thedefinition.}\it}{\medskip}

\newenvironment{eqalph}{\stepcounter{equation}
\setcounter{equationa}{\value{equation}}
\setcounter{equation}{0}

\begin{eqnarray}}{\end{eqnarray}
\setcounter{equation}{\value{equationa}}}

\newcommand{\mar}[1]{}

\begin{document}
\hbox{}

{\parindent=0pt

{\large \bf Lagrangian symmetries and
supersymmetries depending on derivatives. Global analysis}
\bigskip

{\bf Giovanni Giachetta},\footnote{E-mail:
giovanni.giachetta@unicam.it} {\bf Luigi 
Mangiarotti},\footnote{E-mail:
luigi.mangiarotti@unicam.it} 


Department of Mathematics and Informatics, University
of Camerino, 62032 Camerino (MC), Italy

\medskip

{\bf
Gennadi Sardanashvily}\footnote{E-mail:
sard@grav.phys.msu.su}

 Department of Theoretical Physics, Physics Faculty,
Moscow State University, 117234 Moscow, Russia
\bigskip


{\bf Abstract:} 
Generalized symmetries and
supersymmetries depending on derivatives of dynamic variables
are treated in a most general setting.
Studding cohomology of the variational bicomplex, we state the
first variational formula and conservation laws  for Lagrangian
systems on fiber bundles and graded manifolds under generalized
symmetries and supersymmetries of any order. Cohomology of
nilpotent generalized supersymmetries are obtained. 
}

\section{Introduction}

Symmetries of differential equations under transformations
of dynamic variables depending on their derivatives 
have been intensively investigated (see
\cite{and93,ibr,kras,olv} for a survey). 
Following \cite{and93,olv}, we agree to call them the
generalized symmetries in contrast with the classical (point)
ones. In mechanics, conservation laws corresponding to
generalized symmetries are well known
\cite{olv}. In field theory, BRST transformations provide the
most interesting example of generalized supersymmetries
\cite{fat,fulp}. 
 
Generalized symmetries of Lagrangian
systems on a local
coordinate domain have been described
in detail \cite{olv}. We aim to provide the global analysis of 
Lagrangian systems on fiber bundles and graded
manifolds under generalized symmetries and supersymmetries of
any order.

Let us note that an $m$-order differential equation on a
finite-dimensional smooth fiber bundle
$\pi:Y\to X$ is conventionally defined as a closed subbundle of
the
$m$-order jet bundle $J^mY\to X$
of sections of $Y\to X$ \cite{bry,kras}.
Euler--Lagrange equations need not satisfy this condition,
unless an Euler--Lagrange operator is of constant rank.
Therefore, we regard infinitesimal symmetry transformations
as derivations of the graded differential algebra (henceforth
GDA) $\cO^*_\infty$ of exterior
forms on jet manifolds, but not as
manifold maps.  This approach is straightforwardly extended to
Lagrangian systems on graded manifolds. 

We use the first variational formula of the calculus of
variations in order to obtain Lagrangian conservation laws
\cite{book,book00,epr}.  Recall that an $r$-order Lagrangian
on a fiber bundle $Y\to X$ is defined as a horizontal density
$L:J^rY\to
\op\w^nT^*X$,
$n=\di X$, on the $r$-order jet manifold $J^rY$.
If $X=\Bbb R$, we are in the case of non-relativistic
time-dependent mechanics.  A classical symmetry is represented
by a
projectable vector field $u$ on $Y\to X$ seen as
an infinitesimal generator of
a local one-parameter group of bundle automorphisms 
of $Y\to X$. Let $\bL_{J^ru}L$ be the Lie derivative of
$L$ along the jet prolongation $J^ru$ of $u$ onto $J^rY$.
The first
variational formula provides its canonical decomposition 
\mar{g16}\beq
\bL_{J^ru}L=u_V\rfloor\dl L + d_H(h_0(J^{2r-1}u\rfloor\Xi_L)), 
\label{g16}
\eeq
where $\dl L$ is the Euler--Lagrange operator, $\Xi_L$ is
a Lepagean equivalent (e.g., a Poincar\'e--Cartan form) of $L$,
$u_V$ is the vertical part of
$u$,
$d_H$ is the total differential, and $h_0$ is the
horizontal projection (see all the definitions below).  
Let $u$ be a divergence symmetry of $L$, i.e., the Lie
derivative 
$\bL_{J^ru}L$ is a total differential $d_H\si$.
Then the first variational
formula (\ref{g16}) on the  kernel Ker$\,\dl L$ of
the Euler--Lagrange operator $\dl L$ leads to
the weak conservation law
\mar{g19}\beq
0\ap d_H(h_0(J^{2r-1}u\rfloor\Xi_L)-\si). \label{g19}
\eeq
If $u$ is a variational symmetry  of $L$, i.e.,
$\bL_{J^ru}L=0$, the conservation law (\ref{g19})
comes to the conservation law of the Noether current 
$\cJ_u=h_0(J^{2r-1}u\rfloor\Xi_L)$.  

In the case of a
classical symmetry, the first variational formula (\ref{g16})
and the existence of a globally defined symmetry current
$\cJ_u$ issue from the existence of
a global Lepagean equivalent $\Xi_L$ of $L$ \cite{got}. In
order to extend the first variational formula to generalized
symmetries (formula (\ref{g8})) and generalized
supersymmetries (formula (\ref{g107})), we derive it from
the decomposition 
\mar{+421}\beq
dL=\dl L - d_H(\Xi), \label{+421}
\eeq
provided by the global exactness of the
subcomplex of one-contact forms of the variational bicomplex
on fiber bundles and graded manifolds (Propositions
\ref{g93} and \ref{g103}). As a consequence, the
existence of a global finite order Lepagean equivalent of a
graded Lagrangian is stated.

A vector field 
$u$ in the first variational formula (\ref{g16}) represents a
derivation of the $\Bbb R$-ring
$C^\infty(Y)$ of smooth real functions on $Y$.
Accordingly, a $k$-order generalized vector field $\vt$ can be
defined  as a derivation of the $\Bbb R$-ring $C^\infty(Y)$ 
with values into the ring $C^\infty(J^kY)$ of
smooth real functions on the jet manifold $J^kY$.
This definition
recovers the geometric notion of a generalized
vector field as a section of the pull-back bundle
$TY\times J^kY\to J^kY$ in \cite{fat}.
The prolongation $J^r\vt$ of $\vt$  
onto any finite order jet manifold $J^rY$ is that one calls a
generalized symmetry. We give the intrinsic definition
of a generalized symmetry as a derivation $\up$ of the $\Bbb
R$ ring $\cO^0_\infty$ such that the Lie derivative
$\bL_\up$ preserves the ideal of contact forms of the above
mentioned 
GDA $\cO^*_\infty$ (Propositions
\ref{g62} -- \ref{g72}). The key point is that 
the Lie derivative $\bL_{J^r\vt}$ along a generalized symmetry 
$J^r\vt$ sends exterior forms on the
jet manifold $J^r Y$ onto those on the jet
manifold $J^{r+k}Y$. By virtue of the well-known B\"acklund
theorem,
$\bL_{J^r\vt}$ preserves $\cO^*_r$  
iff either $\vt$ is a vector field on $Y$ or $Y\to X$ is a
one-dimensional bundle and $\vt$ is a generalized vector field
at most of first order \cite{ibr}. Thus, considering generalized
symmetries, one deals with Lagrangian systems of unspecified
finite order. 

Infinite order jet formalism provides a
convenient tool for studying these systems
\cite{ander,jmp,kras,book00,epr,tak2}. 
 With the inverse system of finite order jet manifolds
\mar{5.10}\beq
X\op\longleftarrow^\pi Y\op\longleftarrow^{\pi^1_0} J^1Y
\longleftarrow \cdots J^{r-1}Y \op\longleftarrow^{\pi^r_{r-1}}
J^rY\longleftarrow\cdots,
\label{5.10}
\eeq
we have the direct system
\mar{5.7}\beq
\cO^*(X)\op\longrightarrow^{\pi^*} \cO^*(Y) 
\op\longrightarrow^{\pi^1_0{}^*} \cO_1^* \cdots
\op\longrightarrow^{\pi^r_{r-1}{}^*}
 \cO_r^* \longrightarrow\cdots \label{5.7}
\eeq
of GDAs of exterior forms on
these manifolds with respect to 
the pull-back monomorphisms $\pi^r_{r-1}{}^*$. Its direct
limit
 is the above mentioned GDA $\cO_\infty^*$
consisting of all the exterior forms on
finite order jet manifolds modulo the pull-back identification.
The exterior differential on $\cO_\infty^*$ is decomposed into
the sum $d=d_H+d_V$ of the total and the vertical 
differentials. These differentials and
the variational operator
$\dl$ split $\cO_\infty^*$ into the variational
bicomplex (\ref{7}), which provides the algebraic description of
Lagrangian systems on a fiber bundle $Y\to X$.
Restricted to a coordinate domain of 
$Y$, this bicomplex, except the terms $\Bbb R$, is exact.  
One refers to this fact as the algebraic
Poincar\'e lemma  (e.g., \cite{olv}). Recently, we have stated 
cohomology of the variational bicomplex for an arbitrary $Y$
\cite{lmp,jmp,ijmms}. The key point is that this cohomology
provides a topological obstruction to local divergence
symmetries to be the global ones. For instance, if a
generalized symmetry
$\vt$ is a divergence symmetry of a Lagrangian
$L$, the equality
$\dl (\bL_{J^r\vt}L)=0$
holds, but the converse is not true because of the de Rham
cohomology group $H^n(Y)$ of $Y$.

\begin{rem} \label{g200} \mar{g200}
Let us point out the following technical detail repeatedly met
in the sequel. The de Rham cohomology of $\cO_\infty^*$ 
is easily proved to equal the
de Rham cohomology $H^*(Y)$ of $Y$
\cite{ander}. However, one has to enlarge the 
GDA $\cO_\infty^*$
in order to find its $d_H$- and $\dl$-cohomology.
Let $\gO^*_r$ be the sheaf
of germs of exterior forms on the $r$-order jet 
manifold $J^rY$, and let
$\ol\gO^*_r$ be its canonical presheaf. We throughout follow
the sheaf terminology of \cite{hir}. There is the direct  system
of presheaves
\be
\ol\gO^*_X\op\longrightarrow^{\pi^*} \ol\gO^*_0 
\op\longrightarrow^{\pi^1_0{}^*} \ol\gO_1^* \cdots
\op\longrightarrow^{\pi^r_{r-1}{}^*}
 \ol\gO_r^* \longrightarrow\cdots. 
\ee
Its direct limit $\ol\gO^*_\infty$ 
is a presheaf of GDAs on the projective limit
of the inverse system (\ref{5.10}) of jet
manifolds. This projective limit,
called the infinite order jet space  $J^\infty Y$,
is endowed with the weakest
topology such that surjections $\pi^\infty_r:J^\infty Y\to
J^rY$ are continuous. This
topology makes $J^\infty Y$ into a paracompact Fr\'echet
manifold \cite{tak2}. 
Let
$\gQ^*_\infty$ be a sheaf constructed from 
$\ol\gO^*_\infty$. The module 
$\cQ^*_\infty=\G(\gQ^*_\infty)$ of 
sections of $\gQ^*_\infty$ is a GDA such that, given an
element
$\f\in \gQ^*_\infty$ and a point $z\in J^\infty Y$, there 
exist an open
neighbourhood $U$ of $z$ and an
exterior form
$\f^{(k)}$ on some finite order jet manifold $J^kY$ so that
$\f|_U= \f^{(k)}\circ \pi^\infty_k|_U$. In particular, there is
the  monomorphism $\cO^*_\infty
\to\cQ^*_\infty$. The key point is that the paracompact
 space
$J^\infty Y$ admits a partition of unity by elements of
the ring $\cQ^0_\infty$ \cite{tak2} and $Y$ is a strong
deformation retract of $J^\infty Y$
 \cite{ander,jmp}. These facts have enabled one to
obtain 
$d_H$- and $\dl$-cohomology of $\cQ^*_\infty$
\cite{and,ander,tak2}. Recently, we have shown that its
subalgebra 
$\cO^*_\infty$ possesses the same $d_H$- and $\dl$-cohomology as
$\cQ^*_\infty$ \cite{lmp,ijmms}.
\end{rem}

The following two peculiarities of generalized
supersymmetries should be additionally noted. Firstly,
generalized supersymmetries are expressed into jets of odd
variables, and we should define an algebra where they act.
Secondly, generalized supersymmetries can be nilpotent.

Stimulated by BRST theory, we do not focus on particular
geometric models of ghost fields in gauge theories (e.g.,
\cite{bon,sch}), but consider Lagrangian systems of odd
variables in a general setting.  For this purpose, one calls
into play fiber bundles over graded manifolds and
supermanifolds 
\cite{cari,cia,mont}. However, the   
antifield BRST theory on $X=\Bbb
R^n$ \cite{barn95,barn,bran,bran01} involves jets of odd fields
only with respect to space-time coordinates. Therefore, we
describe odd variables on a smooth manifold
$X$ as generating elements of the structure ring of
a graded manifold, whose body is
$X$. By the well-known Batchelor theorem \cite{bart},
such a graded manifold is
isomorphic to the one whose structure sheaf 
$\gA_Q$ is formed by germs of sections of the exterior product 
\mar{g80}\beq
\w Q^*=\Bbb R\op\oplus_X
Q^*\op\oplus_X\op\w^2 Q^*\op\oplus_X\cdots,
\label{g80}
\eeq 
where $Q^*$ is the dual of some real vector
bundle $Q\to X$. In physical models, a vector bundle $Q$ is
usually given from the beginning. Therefore, we restrict our
consideration to graded manifolds
$(X,\gA_Q)$ where the
Batchelor isomorphism holds fixed, i.e., automorphisms of
$(X,\gA_Q)$ are restricted to those induced by bundle
automorphisms of $Q$. This restriction enables us to handle
the structures which are not preserved by general
automorphisms of a graded manifold. We agree to call 
$(X,\gA_Q)$ a simple graded manifold constructed from $Q$.
Accordingly, 
$r$-order jets of odd variables are defined as generating
elements of the structure ring of the simple graded manifold
$(X,\gA_{J^rQ})$ constructed from the jet bundle $J^rQ$ of
$Q$ \cite{book00,mpla}. 
Let $\cC^*_{J^rQ}$ be the bigraded differential algebra
(henceforth BGDA) of graded exterior forms on the graded
manifold
$(X,\gA_{J^rQ})$. Since $\pi^r_{r-1}:J^rQ \to J^{r-1}Q$ is
a linear bundle morphism over $X$,  it yields the
morphism of graded manifolds 
$(X,\gA_{J^rQ})\to (X,\gA_{J^{r-1}Q})$ and the
monomorphism of the BGDAs
$\cC^*_{J^{r-1}Q}\to \cC^*_{J^rQ}$ \cite{bart,book00}. Hence,
there is the direct system of BGDAs 
\mar{g205}\beq
\cC^*_Q\ar^{\pi^{1*}_0} \cC^*_{J^1Q}\ar\cdots 
\cC^*_{J^rQ}\ar^{\pi^{r+1*}_r}\cdots, \label{g205}
\eeq
whose direct limit $\cC^*_\infty$
consists 
of graded exterior forms on graded manifolds
$(X,\gA_{J^rQ})$,
$0\leq r$, modulo the pull-back identification. 

This definition of odd jets 
differs from that of jets of a graded
fiber bundle in \cite{hern}, but reproduces 
the heuristic notion of jets of ghosts
in the above mentioned antifield BRST theory on $\Bbb R^n$.  
Moreover, it enables one to describe odd and
even variables (e.g., classical fields,
ghosts, ghosts-for-ghosts and antifields in BRST theory)
on the same footing. Namely, let  a smooth fiber
bundle $Y\to X$ be affine. Then its de Rham cohomology equals
that of $
X$. Let
$\cP^*_\infty$ be the
$C^\infty(X)$-subalgebra
of the GDA $\cO^*_\infty$ which consists of
exterior forms whose coefficients are polynomial in the fiber
coordinates on $J^\infty Y\to X$. 
Let us consider the product $\cS^*_\infty$ of
graded algebras $\cC_\infty^*$ and $\cP^*_\infty$ over their
common subalgebra $\cO^*(X)$. It is a BGDA. For the sake
of brevity, we continue to call its elements the graded 
forms. 

Similarly to
$\cO^*_\infty$, the BGDA $\cS^*_\infty$ is split into the
graded variational bicomplex, which provides the algebraic
description of Lagrangian systems of even and odd variables
indexed by elements of the fiber bundles $Y$ and $Q$ over a
smooth manifold
$X$. Following a procedure  similar to that in Remark
\ref{g200}, we obtain cohomology of some complexes of the BGDA
$\cS^*_\infty$. These  are the short variational complex of
horizontal (local in the terminology of
\cite{barn,bran}) graded exterior forms (\ref{g111}), the
complex of one-contact graded exterior forms (\ref{g112}) and
the de Rham complex
(\ref{g110}). 
Cohomology of the first and third complexes is proved to equal
the de Rham cohomology $H^*(X)$ of $X$, while the second one is
globally exact (Theorem \ref{g96}). This exactness provides the
decomposition (\ref{+421}) of  a graded Lagrangian and leads to
the first variational formula and conservation law under
generalized supersymmetries (Propositions \ref{g103} and
\ref{g106}). Cohomology
$H^{<n}(X)$ of the short variational complex is the main
ingredient in a computation of the iterated
cohomology of nilpotent generalized supersymmetries.

By analogy with a generalized symmetry, a generalized
supersymmetry $\up$ is defined as a graded derivation of the
$\Bbb R$-ring
$\cS^0_\infty$ such that the Lie derivative $\bL_\up$
preserves the contact ideal of the
BGDA
$\cS^*_\infty$ (Proposition \ref{g231}). 
The BRST transformation
$\up$ (\ref{g130}) in gauge theory on a principal bundle
exemplifies such a generalized supersymmetry. Its peculiarity
is that the Lie derivative $\bL_\up$  of
horizontal graded forms  is nilpotent. This fact 
motivates us to study nilpotent
generalized supersymmetries in a general setting. 

Note that nilpotent generalized supersymmetries are
necessarily odd, i.e., there are no nilpotent generalized
symmetries. The key point is that the Lie derivative
$\bL_\up$ along a generalized supersymmetry and the total
differential
$d_H$ mutually commute. If $\bL_\up$ is nilpotent, 
let us suppose that the
$d_H$-complex $\cS^{0,*}_\infty$ of horizontal graded forms 
is split into a complex of complexes $\{S^{k,m}\}$ with
respect to $\bL_\up$ and
$d_H$. In order to make it into a bicomplex, let us introduce
the nilpotent operator
\mar{g240}\beq
\bs_\up\f=(-1)^{|\f|}\bL_\up\f,
\qquad \f\in S^{0,*}_\infty, \label{g240}
\eeq
such that $d_H\circ\bs_\up+\bs_\up\circ d_H=0$. In the case
of the BRST transformation
$\up$ (\ref{g130}), $\bs_\up$ (\ref{g240}) is the BRST operator.
The bicomplex $S^{*,*}$
is graded by the form degree $0\leq m\leq n$ and an
integer
$k\in\Bbb Z$. Let us consider horizontal graded forms $\f\in
S^{0,*}_\infty$ such that a nilpotent generalized supersymmetry
$\up$ is their divergence symmetry, i.e. $\bL_\up\f=d_H\f$. We
come to the relative and iterated cohomology of
$\bs_\up$ with respect to the total differential $d_H$.
In the antifield BRST theory, relative cohomology is 
known as the local BRST cohomology 
\cite{barn,bran} (see \cite{dub} for the BRST cohomology modulo
the exterior differential $d$).  Relative and iterated
cohomology groups coincide with each other on horizontal
densities, and they naturally characterize graded Lagrangians
$L$, for which $\up$ is a divergence symmetry, modulo the Lie
derivatives $\bL_\up\xi$, $\xi\in S^{0,*}_\infty$, and the
$d_H$-exact graded forms.

We obtain the iterated cohomology $H^{*,m<n}(\bs_\up|d_H)$
(Theorem \ref{g250})
and state the relation between the iterated cohomology
$H^{*,n}(\bs_\up|d_H)$ and the total
$(\bs_\up+d_H)$-cohomology of the bicomplex
$S^{*,*}$ (Theorem \ref{g251}). 
This relation plays a prominent role, e.g., in the antifield
BRST theory \cite{barn,bran}. Note that relative cohomology of
form degree
$m<n$ fails to be related to the total cohomology. For
instance, in Section 9.6 of
\cite{barn}, iterated BRST cohomology in fact is considered.

\section{Lagrangian systems of unspecified finite order on
fiber bundles}

This Section addresses the basic formulae for finite order
Lagrangian systems on a smooth fiber bundle $Y\to X$ in the
framework of infinite order jet formalism. The similar formulae 
for graded Lagrangian systems will be
stated in Section 4. Our main goal is the decomposition
(\ref{+421}).

\begin{rem} \label{g225} \mar{g225}
Smooth manifolds throughout are real, finite-dimensional,
Hausdorff, second-countable (hence, paracompact) and connected. 
\end{rem}

Any bundle coordinate atlas $\{(U_Y;x^\la,y^i)\}$ of $\pi:Y\to
X$ yields the coordinate atlas 
\mar{jet1}\beq
\{((\pi^\infty_0)^{-1}(U_Y); x^\la, y^i_\La)\}, \qquad
{y'}^i_{\la+\La}=\frac{\dr x^\m}{\dr x'^\la}d_\m y'^i_\La,
\qquad
0\leq|\La|,
\label{jet1}
\eeq
of the Fr\'echet manifold $J^\infty Y$, where
$\La=(\la_k...\la_1)$ is a symmetric multi-index 
of length $k$, $\la+\La=(\la\la_k...\la_1)$, 
and  
\mar{5.177}\beq
d_\la = \dr_\la + \op\sum_{|\La|\geq 0}
y^i_{\la+\La}\dr_i^\La, \qquad
d_\La=d_{\la_r}\circ\cdots\circ d_{\la_1}, \quad
\La=(\la_r...\la_1), \label{5.177}
\eeq
are the total derivatives. Hereafter, we 
fix an atlas of $Y$ and, consequently, that of $J^\infty Y$
containing a finite number of charts $(U_Y;x^\la,y^i)$ (their
branches $U_Y$ however need not be domains) \cite{greub}. 

Restricted to a
coordinate chart (\ref{jet1}), elements of the GDA
$\cO^*_\infty$ can be written in a
coordinate form; horizontal forms 
$\{dx^\la\}$ and contact one-forms
$\{\th^i_\La=dy^i_\La -y^i_{\la+\La}dx^\la\}$ make up a local
basis for the $\cO^0_\infty$-algebra
$\cO^*_\infty$. 
There is the canonical decomposition 
$\cO^*_\infty =\oplus\cO^{k,m}_\infty$
of this algebra into $\cO^0_\infty$-modules $\cO^{k,m}_\infty$
of $k$-contact and $m$-horizontal forms
together with the corresponding
projections $h_k:\cO^*_\infty\to \cO^{k,*}_\infty$ and
$h^m:\cO^*_\infty\to \cO^{*,m}_\infty$.
Accordingly, the
exterior differential on $\cO_\infty^*$ is split
into the sum $d=d_H+d_V$ of the total and vertical
differentials 
\be
&& d_H\circ h_k=h_k\circ d\circ h_k, \qquad d_H\circ
h_0=h_0\circ d, \qquad d_H(\f)= dx^\la\w d_\la(\f), \\ 
&& d_V \circ h^m=h^m\circ d\circ h^m, \qquad
d_V(\f)=\th^i_\La \w \dr^\La_i\f, \qquad \f\in\cO^*_\infty.
\ee
One also introduces the $\Bbb R$-module
projector 
\mar{r12}\beq
\vr=\op\sum_{k>0} \frac1k\ol\vr\circ h_k\circ h^n,
\qquad \ol\vr(\f)= \op\sum_{|\La|\geq 0}
(-1)^{\nm\La}\th^i\w [d_\La(\dr^\La_i\rfloor\f)], 
\qquad \f\in \cO^{>0,n}_\infty, \label{r12}
\eeq
of $\cO^*_\infty$ such that
$\vr\circ d_H=0$ and the nilpotent variational operator
$\dl=\vr\circ d$ on $\cO^{*,n}_\infty$. Put
$E_k=\vr(\cO^{k,n}_\infty)$. As a consequence, the GDA 
$\cO^*_\infty$ is split into the above mentioned variational
bicomplex
\mar{7}\beq
\begin{array}{ccccrlcrlccrlccrlcrl}
 & &  &  & & \vdots & & & \vdots  & & & 
&\vdots  & & & &
\vdots & &   & \vdots \\ 
& & & & _{d_V} & \put(0,-7){\vector(0,1){14}} & & _{d_V} &
\put(0,-7){\vector(0,1){14}} & &  & _{d_V} &
\put(0,-7){\vector(0,1){14}} & & &  _{d_V} &
\put(0,-7){\vector(0,1){14}}& & _{-\dl} & \put(0,-7){\vector(0,1){14}} \\ 
 &  & 0 & \to & &\cO^{1,0}_\infty &\ar^{d_H} & &
\cO^{1,1}_\infty &
\ar^{d_H} &\cdots  & & \cO^{1,m}_\infty &\ar^{d_H} &\cdots & &
\cO^{1,n}_\infty &\ar^\vr &  & E_1\to  0\\  
& & & & _{d_V} &\put(0,-7){\vector(0,1){14}} & & _{d_V} &
\put(0,-7){\vector(0,1){14}} & & &  _{d_V}
 & \put(0,-7){\vector(0,1){14}} & &  & _{d_V} & \put(0,-7){\vector(0,1){14}}
 & & _{-\dl} & \put(0,-7){\vector(0,1){14}} \\
0 & \to & \Bbb R & \to & & \cO^0_\infty &\ar^{d_H} & &
\cO^{0,1}_\infty &
\ar^{d_H} &\cdots  & &
\cO^{0,m}_\infty & \ar^{d_H} & \cdots & &
\cO^{0,n}_\infty & \equiv &  & \cO^{0,n}_\infty \\
& & & & _{\pi^{\infty*}}& \put(0,-7){\vector(0,1){14}} & & _{\pi^{\infty*}} &
\put(0,-7){\vector(0,1){14}} & & &  _{\pi^{\infty*}}
 & \put(0,-7){\vector(0,1){14}} & &  & _{\pi^{\infty*}} &
\put(0,-7){\vector(0,1){14}} & &  & \\
0 & \to & \Bbb R & \to & & \cO^0(X) &\ar^d & & \cO^1(X) &
\ar^d &\cdots  & &
\cO^m(X) & \ar^d & \cdots & &
\cO^n(X) & \ar^d & 0 &  \\
& & & & &\put(0,-5){\vector(0,1){10}} & & &
\put(0,-5){\vector(0,1){10}} & & & 
 & \put(0,-5){\vector(0,1){10}} & & &   &
\put(0,-5){\vector(0,1){10}} & &  & \\
& & & & &0 & &  & 0 & & & & 0 & & & & 0 & &  & 
\end{array}
\label{7}
\eeq

The second row from the bottom and the last column of this
bicomplex assemble into the variational complex
\mar{b317}\beq
0\to\Bbb R\to \cO^0_\infty
\ar^{d_H}\cO^{0,1}_\infty\cdots  
\op\longrightarrow^{d_H} 
\cO^{0,n}_\infty  \op\longrightarrow^\dl E_1 
\op\longrightarrow^\dl 
E_2 \ar \cdots\,. \label{b317} 
\eeq
One can think of its elements
\be
L=\cL\om\in \cO^{0,n}_\infty, 
\qquad \dl L=\op\sum_{|\La|\geq
0}(-1)^{|\La|}d_\La(\dr^\La_i \cL)\th^i\w\om\in E_1, \qquad
 \om=dx^1\w\cdots\w dx^n,
\ee
as being
a finite order Lagrangian and its
Euler--Lagrange operator.

\begin{theo} \label{g90} \mar{g90} 
Cohomology of the variational complex 
(\ref{b317}) is isomorphic to the de Rham cohomology of a
fiber bundle
$Y$ \cite{lmp,ijmms}. 
\end{theo}

\noindent {\it Outline of proof.}  We have the complex of
sheaves of
$\cQ^0_\infty$-modules 
\mar{g91}\beq
0\to\Bbb R\to \gQ^0_\infty
\ar^{d_H}\gQ^{0,1}_\infty\cdots  
\op\longrightarrow^{d_H} 
\gQ^{0,n}_\infty  \op\longrightarrow^\dl \gE_1 
\op\longrightarrow^\dl 
\gE_2 \longrightarrow \cdots, \qquad
\gE_k=\vr(\gO^{k,n}_\infty), 
\label{g91} 
\eeq
on $J^\infty Y$ and the complex of their
structure modules
\mar{g92}\beq
0\to\Bbb R\to \cQ^0_\infty
\ar^{d_H}\cQ^{0,1}_\infty\cdots  
\op\longrightarrow^{d_H} 
\cQ^{0,n}_\infty  \op\longrightarrow^\dl \cE_1 
\op\longrightarrow^\dl 
\cE_2 \longrightarrow \cdots\,.  \label{g92} 
\eeq
Since the paracompact space $J^\infty Y$ admits 
a partition of unity by elements of the ring $\cQ^0_\infty$,  
the sheaves $\gQ^{0,k}_\infty$ in the complex (\ref{g91}) are
fine. The sheaves $\gE_k$ are also proved to be fine
\cite{lmp,ijmms}. Consequently, all sheaves, except $\Bbb
R$, in the complex (\ref{b317}) are acyclic. Then,
by virtue of the above mentioned algebraic Poincar\'e lemma,
the  complex (\ref{g91}) is a resolution of the constant sheaf
$\Bbb R$ on $J^\infty Y$. In accordance with the abstract de
Rham theorem \cite{hir}, cohomology of
the complex (\ref{g92}) equals the cohomology 
of $J^\infty Y$ with coefficients in $\Bbb R$. The
latter, in turn, is isomorphic to  the de Rham cohomology
of $Y$, which is a strong deformation retract of $J^\infty Y$.
Finally, the $d_H$- and $\dl$-cohomology of
$\cQ^*_\infty$ is proved to equal that of its
subalgebra 
$\cO^*_\infty$ \cite{lmp,ijmms}. $\Box$ \medskip

\begin{cor} \label{g212} \mar{g212}
Every $d_H$-closed form  $\f\in 
\cO^{0,m<n}$ is the sum 
\mar{g213}\beq
\f=h_0\vf +d_H\xi, \qquad \xi\in
\cO^{0,m-1}_\infty, \label{g212}
\eeq
where $\vf$ is a closed $m$-form on $Y$. 
Every $\dl$-closed form (a variationally trivial Lagrangian)
$L\in\cO^{0,n}$ is the sum
\mar{t42}\beq
 L=h_0\varphi + d_H\xi,  \qquad \xi\in
\cO^{0,n-1}_\infty,
\label{t42}
\eeq
where $\varphi$ is a closed $n$-form on $Y$.
\end{cor} 

\begin{rem}
The formulae (\ref{g212}) -- (\ref{t42})  have been stated in
\cite{and} by computing cohomology of the fixed order
variational sequence, but the proof of the local exactness of
this sequence requires rather sophisticated {\it ad hoc}
techniques. The proof of Theorem \ref{g96} below on cohomology
of the graded variational bicomplex follows the above proof of 
Theorem \ref{g90}.
\end{rem} 

\begin{prop} \label{g93} \mar{g93}
For any Lagrangian $L\in \cO^{0,n}_\infty$, there is the
decomposition (\ref{+421}), where $\Xi\in \cO^{1,n-1}_\infty$.
\end{prop}

\begin{proof} Let us consider 
the third row  
\mar{g201}\beq
0\to \cO^{1,0}_\infty\ar^{d_H} \cO^{1,1}_\infty
\cdots
\ar^{d_H}\cO^{1,n}_\infty\ar^\vr E_1\to 0
\label{g201}
\eeq
from the bottom of the variational bicomplex (\ref{7}).
 Similarly to the proof of
Theorem
\ref{g90}, one can show that the
complex  (\ref{g201}) is exact. Its exactness at the term
$\cO^{1,n}_\infty$ relative to the 
projector $\vr$ provides 
the $\Bbb R$-module decomposition
\be
\cO^{1,n}_\infty=E_1\oplus d_H(\cO^{1,n-1}_\infty).
\ee
It leads to the splitting
(\ref{+421}) of
$dL\in\cO^{1,n}_\infty$. The form $\Xi$ in this
splitting is not uniquely defined. It reads
\mar{g43}\ben
&& \Xi=\op\sum_{s=0}F^{\la\nu_s\ldots\nu_1}_i
\th^i_{\nu_s\ldots\nu_1}\w\om_\la,\qquad
\om_\la=\dr_\la\rfloor\om, \label{g43}\\
&& F_i^{\nu_k\ldots\nu_1}=
\dr_i^{\nu_k\ldots\nu_1}\cL-d_\la F_i^{\la\nu_k\ldots\nu_1}
+h_i^{\nu_k\ldots\nu_1},  \nonumber
\een
where local functions $h\in\cO^0_\infty$ obey the relations
$h^\nu_i=0$,
$h_i^{(\nu_k\nu_{k-1})\ldots\nu_1}=0$. It follows that
$\Xi_L=\Xi+L$ is a Lepagean equivalent, e.g., a
Poincar\'e--Cartan form of a finite order Lagrangian $L$
\cite{got}.
\end{proof}

\section{Generalized Lagrangian symmetries}

A derivation $\up\in\gd\cO^0_\infty$ of the
$\Bbb R$-ring $\cO^0_\infty$
is said to be a
generalized symmetry if the Lie derivative $\bL_\up$, being a
derivation of the GDA $\cO^*_\infty$, preserves its ideal of
contact forms.
Forthcoming Propositions \ref{g62} -- \ref{g72} confirm
the contentedness of this definition.

\begin{prop} \label{g62} \mar{g62}
The derivation module $\gd\cO^0_\infty$ is isomorphic to the 
$\cO^0_\infty$-dual $(\cO^1_\infty)^*$ of the module of
one-forms $\cO^1_\infty$.
\end{prop}

\begin{proof}
At first, let us show that $\cO^*_\infty$ is generated by
elements $df$,
$f\in \cO^0_\infty$. It suffices to justify that any
element of $\cO^1_\infty$ is a finite $\cO^0_\infty$-linear
combination of elements $df$,
$f\in \cO^0_\infty$. Indeed, every
$\f\in\cO^1_\infty$ is an exterior form on some finite order
jet manifold $J^rY$.  By
virtue of the Serre--Swan theorem extended to vector
bundles over non-compact manifolds \cite{ren,ss}, the
$C^\infty(J^rY)$-module $\cO^*_r$ of one-forms on $J^rY$ is
a projective module of finite rank, i.e., $\f$ is
represented by a finite $C^\infty(J^rY)$-linear combination of
elements $df$,
$f\in C^\infty(J^rY)\subset
\cO^0_\infty$.  Any element $\Phi\in (\cO^1_\infty)^*$ yields a
derivation
$\up_\Phi(f)=\Phi(df)$ of the $\Bbb R$-ring $\cO^0_\infty$.
Since the module
$\cO^1_\infty$ is generated by elements $df$, $f\in
\cO^0_\infty$, different elements of $(\cO^1_\infty)^*$
provide different derivations of $\cO^0_\infty$, i.e., there
is a monomorphism  $(\cO^1_\infty)^*\to \gd\cO^0_\infty$. 
By
the same formula, any derivation $\up\in \gd\cO^0_\infty$
sends $df\mapsto \up(f)$ and, since $\cO^0_\infty$ is generated
by elements $df$, it defines a morphism 
$\Phi_\up:\cO^1_\infty\to
\cO^0_\infty$. Moreover, different derivations $\up$ provide
different morphisms $\Phi_\up$. Thus, we have a monomorphism
and, consequently, an isomorphism $\gd\cO^0_\infty\to 
(\cO^1_\infty)^*$. 
\end{proof}

\begin{rem}
As follows from Proposition \ref{g62}, the de Rham complex
of the GDA
$\cO^*_\infty$ is both the Chevalley--Eilenberg complex of the
Lie algebra
$\gd\cO^0_\infty$ with coefficients in $\cO^0_\infty$ and the
universal differential calculus over the $\Bbb R$-ring
$\cO^0_\infty$.
\end{rem}

\begin{prop} \label{g60} \mar{g60}
Relative to an atlas
(\ref{jet1}), a derivation $\up\in\gd\cO^0_\infty$ is
given by the expression
\mar{g3}\beq
\up=\up^\la \dr_\la + \up^i\dr_i +
\op\sum_{|\La|>0}\up^i_\La
\dr^\La_i, \label{g3}  
\eeq
where $\up^\la$, $\up^i$, $\up^i_\La$ are local
smooth functions of finite jet order obeying the transformation
law
\mar{g71}\beq
\up'^\la=\frac{\dr x'^\la}{\dr x^\m}\up^\m, \qquad
\up'^i=\frac{\dr y'^i}{\dr y^j}\up^j + \frac{\dr y'^i}{\dr
x^\m}\up^\m, \qquad 
\up'^i_\La=\op\sum_{|\Si|\leq|\La|}\frac{\dr y'^i_\La}{\dr
y^j_\Si}\up^j_\Si +
\frac{\dr y'^i_\La}{\dr x^\m}\up^\m. \label{g71} 
\eeq
\end{prop}

\begin{proof}
Restricted to a coordinate chart (\ref{jet1}), 
$\cO^1_\infty$ is a free $\cO^0_\infty$-module
generated by the
exterior forms $dx^\la$, $\th^i_\La$. Then
$\gd\cO^0_\infty=(\cO^1_\infty)^*$ restricted to this chart 
consists of elements (\ref{g3}), where $\dr_\la$, $\dr^\La_i$
are the duals of $dx^\la$, $\th^i_\La$. The
transformation rule (\ref{g71}) results from the transition
functions (\ref{jet1}). Since the atlas (\ref{jet1}) is finite,
a derivation $\up$ (\ref{g3}) preserves $\cO^*_\infty$.
\end{proof}

The interior product $\up\rfloor\f$ 
and the Lie derivative $\bL_\up\f$,
$\f\in\cO^*_\infty$, obey  
the standard formulae. 

\begin{prop} \label{g72} \mar{g72}
A derivation $\up$ (\ref{g3}) is a generalized
symmetry iff 
\mar{g4}\beq
\up^i_\La=d_\La(\up^i-y^i_\m\up^\m)+y^i_{\m+\La}\up^\m, \qquad
0<|\La|.
\label{g4}
\eeq
\end{prop}

\begin{proof}
The expression (\ref{g4}) results from a direct
computation  similarly to the first part
of the above mentioned B\"acklund theorem. Then one can 
justify that local functions (\ref{g4}) fulfill the 
transformation law (\ref{g71}).
\end{proof}

Any generalized symmetry admits the horizontal
splitting
\mar{g5}\beq
\up=\up_H +\up_V=\up^\la d_\la + (\vt^i\dr_i +
\op\sum_{|\La|>0} d_\La \vt^i\dr_i^\La), \qquad \vt^i=
\up^i-y^i_\m\up^\m,
\label{g5}
\eeq
 relative to the canonical connection
$\nabla=dx^\la\ot d_\la$ on the $C^\infty(X)$-ring
$\cO^0_\infty$ \cite{book00}.
For instance,  let $\tau$ be a vector field on $X$. Then
the derivation 
$\tau\rfloor (d_H f)$, $f\in \cO_\infty^0$, of $\cO_\infty^0$
is a horizontal generalized symmetry
$\up=\tau^\m d_\m$.
A direct computation shows that any
vertical generalized symmetry
$\up=\up_V$ obeys the relations
\mar{g6,'}\ben
&& \up\rfloor d_H\f=-d_H(\up\rfloor\f), \label{g6}\\
&& \bL_\up(d_H\f)=d_H(\bL_\up\f), \qquad \f\in\cO^*_\infty.
\label{g6'}
\een 

\begin{prop}  \label{g75} \mar{g75}
Given a  Lagrangian $L\in\cO^{0,n}_\infty$, its 
Lie derivative $\bL_\up L$
along a generalized
symmetry
$\up$ (\ref{g5}) fulfils the first variational formula
\mar{g8}\beq
\bL_\up L= \up_V\rfloor\dl L +d_H(h_0(\up\rfloor\Xi_L)) 
+\cL d_V (\up_H\rfloor\om), \label{g8}
\eeq
where $\Xi_L$ is a Lepagean equivalent, e.g.,
a Poincar\'e--Cartan form of $L$.
\end{prop}

\begin{proof}
The formula (\ref{g8}) comes from the splitting
(\ref{+421}) and the relation (\ref{g6}) as follows:  
\mar{g7}\ben
&& \bL_\up L=\up\rfloor dL + d(\up\rfloor L)
=\up_V\rfloor dL + d_H(\up_H\rfloor L) +\cL d_V
(\up_H\rfloor\om)= \label{g7} \\
&& \qquad  \up_V\rfloor\dl L -\up_V\rfloor
d_H\Xi + d_H(\up_H\rfloor L) +\cL d_V (\up_H\rfloor\om)
= \nonumber \\
&& \qquad \up_V\rfloor\dl L +d_H(\up_V\rfloor\Xi +
\up_H\rfloor L) +\cL d_V (\up_H\rfloor\om), \qquad 
\nonumber
\een
where we put $\Xi_L=\Xi+L$.
\end{proof}

In comparison with the first variational formula (\ref{g16})
for classical symmetries, the right-hand side of the
first variational formula (\ref{g8}) contains an additional
contact term which vanishes if a generalized symmetry
$\up$ is projected onto $X$, i.e.,  its components $\up^\la$
depend only on coordinates on $X$. 

Let $\up$ be a divergence
symmetry of $L$, i.e., $\bL_\up L=d_H\si$, $\si\in
\cO^{0,n-1}_\infty$.
By virtue of the expression (\ref{g7}), this condition implies
that a generalized symmetry $\up$ is projected onto $X$.
Then the first variational formula (\ref{g8}) takes the form 
\mar{g11}\beq
d_H\si= \up_V\rfloor\dl L +d_H(h_0(\up\rfloor\Xi_L)).
\label{g11}
\eeq
Restricted to Ker$\,\dl L$, it leads to the weak
conservation law
\mar{g32}\beq
0\ap d_H(h_0(\up\rfloor\Xi_L)-\si). \label{g32}
\eeq

A glance at the expression (\ref{g7}) shows
that a generalized symmetry $\up$ (\ref{g5}), projected onto
$X$, is a divergence symmetry of a Lagrangian $L$ iff its
vertical part $\up_V$ is so. Moreover, $\up$ and $\up_V$ lead to
the same conservation law (\ref{g32}). Thus, we can restrict
our consideration to vertical divergence symmetries $\up$. In
this case, the conservation law (\ref{g32}) takes the form
\be
0\ap d_H(\up\rfloor\Xi_L-\si),
\ee
where $\cJ_\up=\up\rfloor\Xi_L$ is the Noether current along
 $\up$. 

It should be noted that a generalized
symmetry is almost never a variational symmetry of a Lagrangian.
Let us obtain the characteristic equation for
divergence symmetries of a Lagrangian $L$. Let
$\up$ be a vertical generalized symmetry. Then the Lie
derivative $\bL_\up L$ (\ref{g7}) is a horizontal density. Let
us require that it is a $\dl$-closed form, i.e., 
$\dl(\bL_\up L)=0$.
In accordance with the equality (\ref{t42}), this condition is
fulfilled iff 
\mar{g95}\beq
\bL_\up L=h_0\vf +d_H\si, \label{g95}
\eeq
where $\vf$ is a closed $n$-form on $Y$, i.e., $\up$  at least
locally is a divergence symmetry of $L$. It is readily
observed that the topological obstruction $h_0\f$ (\ref{g95})
for
$\up$ to be a global divergence symmetry is at most of first
order. If
$Y\to X$ is an affine bundle, its de Rham cohomology equals
that of
$X$ and, consequently, the topological obstruction
$h_0\vf=\vf$  (\ref{g95}) reduces to a non-exact $n$-form on
$X$. Recall that, by virtue of the master identity
\be
\bL_{J^{2r}u}\dl L=\dl(\bL_{J^ru}L),
\ee  
any classical divergence symmetry of a
Lagrangian
is also a symmetry of its Euler--Lagrange operator.
However, this equality is not true for generalized
symmetries \cite{olv}. It comes to the relation  
\be
\dl(\bL_\up L)=\bL_\up\dl L + \op\sum_{|\La|>0}
(-1)^{|\La|}d_\La
(\dr^\La_k\up^i\dl_i\cL dy^k)\w\om. 
\ee

\section{Graded Lagrangian systems}

Let $(X,\gA_Q)$ be the simple graded manifold constructed from
a vector bundle $Q\to X$ of fiber dimension $m$.
Its structure ring
$\cA_Q$ of sections of $\gA_Q$ consists of sections of the
exterior bundle (\ref{g80}) called graded functions. Given
bundle coordinates
$(x^\la,q^a)$ on
$Q$ with transition functions
$q'^a=\rho^a_b q^b$, let
$\{c^a\}$ be the corresponding fiber bases for
$Q^*\to X$, together with transition functions
$c'^a=\rho^a_bc^b$. Then $(x^\la, c^a)$ is called the local
basis for the graded manifold $(X,\gA_Q)$ \cite{bart,book00}.
With respect to this basis, graded functions read 
\be
f=\op\sum_{k=0}^m \frac1{k!}f_{a_1\ldots
a_k}c^{a_1}\cdots c^{a_k}, 
\ee
where $f_{a_1\cdots
a_k}$ are local smooth real functions on $X$.

Given a graded manifold $(X,\gA_Q)$, by the sheaf $\gd\gA_Q$ of
graded derivations of $\gA_Q$ is meant a subsheaf of
endomorphisms of the structure sheaf $\gA_Q$ such that any
section $u$ of $\gd\gA_Q$ over an open subset $U\subset X$ is a
graded derivation of the graded ring $\cA_Q(U)$ of graded
functions on $U$, i.e.,
\be
u(ff')=u(f)f'+(-1)^{[u][f]}fu (f'), 
\qquad f,f'\in \cA_Q(U),
\ee
for homogeneous elements $u\in\gd\gA_Q(U)$ and $f,f'\in
\cA_Q(U)$, where $[.]$ denotes the Grassmann parity. One can
show that sections of $\gd\gA_Q$ over $U$ exhaust all graded
derivations of the graded ring $\cA_Q(U)$
\cite{bart}. Let $\gd\cA_Q$ be the Lie
superalgebra of graded derivations
of the $\Bbb R$-ring $\cA_Q$. 
Its elements are
called graded vector fields on $(X,\gA_Q)$. Due to the
canonical splitting
$VQ= Q\times Q$, the vertical tangent bundle 
$VQ\to Q$ of $Q\to X$ can be provided with the fiber bases 
$\{\dr_a\}$ which is the dual of 
$\{c^a\}$. Then
a graded vector field takes the local form
$u= u^\la\dr_\la + u^a\dr_a$,
where $u^\la, u^a$ are local graded functions. It
acts on
$\cA_Q$ by the rule
\mar{cmp50'}\beq
u(f_{a\ldots b}c^a\cdots c^b)=u^\la\dr_\la(f_{a\ldots b})c^a\cdots c^b +u^d
f_{a\ldots b}\dr_d\rfloor (c^a\cdots c^b). \label{cmp50'}
\eeq
This rule implies the corresponding transformation law 
\be
u'^\la =u^\la, \qquad u'^a=\rho^a_ju^j +
u^\la\dr_\la(\rho^a_j)c^j. 
\ee
Then one can show \cite{book00,ijmp} that graded vector fields 
on a simple graded manifold can be represented
by sections of the vector bundle
$\cV_Q\to X$ which is locally isomorphic to the vector bundle
\be
\cV_Q|_U\approx\w Q^*\op\ot_X(Q\op\oplus_X TX)|_U,
\ee
and is equipped with the bundle coordinates $(x^\la_{a_1\ldots
a_k},v^i_{b_1\ldots b_k})$,
$k=0,\ldots,m$, together with the transition functions
\be
&& x'^\la_{i_1\ldots i_k}=\rho^{-1}{}_{i_1}^{a_1}\cdots
\rho^{-1}{}_{i_k}^{a_k} x^\la_{a_1\ldots a_k}, \\
&& v'^i_{j_1\ldots j_k}=\rho^{-1}{}_{j_1}^{b_1}\cdots
\rho^{-1}{}_{j_k}^{b_k}\left[\rho^i_jv^j_{b_1\ldots b_k}+ \frac{k!}{(k-1)!} 
x^\la_{b_1\ldots b_{k-1}}\dr_\la\rho^i_{b_k}\right].
\ee

Using this fact, one can introduce graded exterior forms on the
graded manifold $(X,\gA_Q)$ as sections of the 
exterior bundle $\op\w^k\cV^*_Q$, where 
$\cV^*_Q\to  X$ is the pointwise $\w Q^*$-dual of $\cV_Q$.
Relative to the dual bases $\{dx^\la\}$ for $T^*X$ and
$\{dc^b\}$ for $Q^*$, sections of $\cV^*_Q\to  X$ (graded
one-forms) read 
\be
\f=\f_\la dx^\la + \f_adc^a,\qquad \f'_a=\rho^{-1}{}_a^b\f_b,
\qquad
\f'_\la=\f_\la +\rho^{-1}{}_a^b\dr_\la(\rho^a_j)\f_bc^j.
\ee
The duality morphism is given by the interior product 
\be
u\rfloor \f=u^\la\f_\la + (-1)^{[\f_a]}u^a\f_a. 
\ee
Graded exterior forms constitute the BGDA $\cC^*_Q$ with
respect to the graded exterior product $\w$ and the even
exterior differential $d$. Recall the standard formulae
\be
&& \f\w\si =(-1)^{\nm\f\nm\si +[\f][\si]}\si\w \f, \qquad
d(\f\w\si)= (d\f)\w\si +(-1)^{\nm\f}\f\w(d\si), \\
&& u\rfloor(\f\w\si)=(u\rfloor \f)\w\si
+(-1)^{|\f|+[\f][u]}\f\w(u\rfloor\si), \\
&& \bL_u\f=u\rfloor d\f+ d(u\rfloor\f), \qquad
\bL_u(\f\w\si)=\bL_u(\f)\w\si +(-1)^{[u][\f]}\f\w\bL_u(\si).
\ee

Since the jet bundle $J^rQ\to X$ of a vector bundle $Q\to X$
is a vector bundle, let us consider the simple graded manifold
$(X,\gA_{J^rQ})$ constructed from 
$J^rQ\to X$. Its local basis is $\{x^\la,c^a_\La\}$, 
$0\leq |\La|\leq r$,
together with the transition functions
\mar{+471}\beq
c'^a_{\la +\La}=d_\la(\rho^a_j c^j_\La),
\qquad 
d_\la=\dr_\la + \op\sum_{|\La|<r}c^a_{\la+\La}
\dr_a^\La, \label{+471}
\eeq 
where $\dr_a^\La$ are the duals of $c^a_\La$.  
Let $\cC^*_{J^rQ}$ be the above mentioned BGDA of graded
exterior forms on the graded manifold $(X,\gA_{J^rQ})$.
It is locally a free
$C^\infty(X)$-algebra finitely generated by the elements 
$(1, c^a_\La, dx^\la,\th^a_\La=dc^a_\La -c^a_{\la
+\La}dx^\la)$, $0\leq|\La|\leq r$. The direct limit
$\cC^*_\infty$ of the direct system (\ref{g205}) inherits
the BGDA operations which commute with the monomorphisms
$\pi^{r+1*}_r$. It is locally a free
$C^\infty(X)$-algebra countably generated by the elements 
$(1, c^a_\La, dx^\la,\th^a_\La)$, $0\leq|\La|$.

It should be emphasized that, in contrast with the GDA
$\cO^*_\infty$, the BGDA $\cC^*_\infty$ consists
of sections of sheaves on $X$. In order to regard these
algebras  on the same footing, let us consider the open
surjection
$\pi^\infty: J^\infty\to X$ and the direct image 
$\pi^*_\infty\gQ^*_\infty$ on $X$ of the sheaf $\gQ^*_\infty$ 
of exterior forms on $J^\infty Y$. Its stalk at a point $x\in
X$ consists of the equivalence classes of sections of the sheaf
$\gQ^*_\infty$ which coincide on the inverse images
$(\pi^\infty)^{-1}(U_x)$ of open neighbourhoods $U_x$ of $X$.
Since $(\pi^\infty)^{-1}(U_x)$ is the infinite order jet space
of sections of the fiber bundle $\pi^{-1}(U_x)\to X$, 
every point
$x\in X$ has a base of open neighbourhoods
$\{U_x\}$ such that the sheaves $\gQ^{*,*}_\infty$ of
$\cQ^0_\infty$-modules and the sheaves $\gE_k$ in 
Theorem (\ref{g90}) are acyclic on the inverse
images
$(\pi^\infty)^{-1}(U_x)$ of these neighbourhoods. Then, in
accordance with the Leray theorem \cite{god}, cohomology of
$J^\infty Y$ with coefficients in the sheaves $\gQ^{*,*}_\infty$
and $\gE_k$ is isomorphic to that of $X$ with coefficients in
their direct images $\pi^*_\infty\gQ^{*,*}_\infty$ and
$\pi^*_\infty\gE_k$, i.e., the sheaves 
$\pi^*_\infty\gQ^{*,*}_\infty$ and $\pi^*_\infty\gE_k$ on $X$
are acyclic. Hereafter, let
$Y\to X$ be an affine bundle. Then $X$ is a strong
deformation retract of $J^\infty Y$. In this case, the  inverse
images 
$(\pi^\infty)^{-1}(U_x)$ of contractible neighbourhoods $U_x$
are contractible and $\pi^\infty_*\Bbb R=\Bbb R$. Then, by 
virtue of the algebraic Poincar\'e lemma, the variational
bicomplex
$\gQ^*_\infty$ of sheaves on 
$(\pi^\infty)^{-1}(U_x)$, except the terms $\Bbb R$, is exact,
and the variational bicomplex
$\pi^\infty_*\gQ^*_\infty$ of sheaves on $X$ is so. There is
the $\Bbb R$-algebra isomorphism of the GDA of sections of the
sheaf $\pi^\infty_*\gQ^*_\infty$ on $X$ to the GDA
$\cQ^*_\infty$. Thus, the GDA
$\cQ^*_\infty$ and its subalgebra
$\cO^*_\infty$ can be regarded as algebras  of sections of
sheaves on $X$, and they keep their
$d$-, $d_H$- and $\dl$-cohomology expressed into the de Rham
cohomology $H^*(X)=H^*(Y)$ of $X$ \cite{lmp}. 

Let us restrict our consideration to the above mentioned
polynomial subalgebra $\cP^*_\infty$ of $\cO^*_\infty$.
Let us consider the graded product $\cS^*_\infty=
\cC_\infty^*\w\cP^*_\infty$ of
graded algebras $\cC_\infty^*$ and $\cP^*_\infty$ over their
common graded subalgebra $\cO^*(X)$. It consists of the
elements 
\be
&& \psi\ot\f, \qquad (\psi\w\si)\ot\f=\psi\ot(\si\w\f),\qquad 
\psi\in \cC^*_\infty, \quad \f\in \cP^*_\infty,\quad \si\in
\cO^*_\infty,\\
&&\f\ot\psi, \qquad (\f\w\si)\ot\psi=\f\ot(\si\w\psi),\qquad 
\psi\in \cC^*_\infty, \quad \f\in \cP^*_\infty,\quad \si\in
\cO^*_\infty,
\ee
of the tensor products $\cC_\infty^*\ot
\cP^*_\infty$ and
$\cP_\infty^*\ot \cC^*_\infty$ of the
$\cO^*_\infty$-modules 
$\cC_\infty^*$ and $\cP^*_\infty$ which are subject to the
commutation relation
\be
\psi\ot\f=(-1)^{|\psi||\f|}\f\ot\psi
\ee
and the multiplication
\be
(\psi\ot\f)\w(\psi'\ot\f')=(-1)^{|\psi'||\f|}(\psi\w\psi')\ot
(\f\w\f'),
\ee
written for homogeneous elements of graded algebras 
$\cC_\infty^*$ and $\cP^*_\infty$. Introducing the 
notation 
\be
\psi\ot 1=1\ot\psi=\psi, \qquad 1\ot\f=\f\ot 1=\f,
\qquad \psi\ot\f=(\psi\ot 1)\w(1\ot\f)=\psi\w\f,
\ee
one can think of $\cS^*_\infty$ as being a bigraded algebra
generated by elements of $\cC_\infty^*$ and $\cP^*_\infty$ 
and provided with the total form degree $|\psi\w\f|=|\psi|+|\f|$
and the total Grassmann parity $[\psi\w\f]=[\psi]$.
For instance, elements
of the ring
$S^0_\infty$ are polynomials of
$c^a_\La$ and $y^i_\La$ with coefficients in $C^\infty(X)$. 
The
sum of exterior differentials on
$\cC_\infty^*$ and $\cP^*_\infty$ makes $\cS^*_\infty$ into a
BGDA with the standard rules
\be
\vf\w\vf' =(-1)^{|\vf||\vf'| +[\vf][\vf']}\vf'\w \vf, \qquad
d(\vf\w\vf')= d\vf\w\vf' +(-1)^{|\vf|}\vf\w d\vf'.
\ee
It is locally generated by the elements 
$(1, c^a_\La, y^i_\La, dx^\la,\th^a_\La,\th^i_\La)$, $|\La|\geq
0$.

\begin{rem} \label{g207} \mar{g207} 
If $Y\to X$ is a vector bundle, one can get the BGDA 
$\cS^*_\infty=
\cC_\infty^*\w\cP^*_\infty$ in a different way.
 Let us consider the Whitney sum 
$S=Q\oplus Y$  of
vector bundles $Q\to X$ and $Y\to X$ regarded 
as a bundle of graded vector spaces $Q_x\oplus Y_x$,
$x\in X$. Let us define the quotient
$\ol S^k$ of the tensor product 
\be
S^k=\Bbb R\op\oplus_X
S^*\op\oplus_X\op\ot^2 S^*\op\oplus_X\cdots \op\oplus_X \op\ot^k
S^* 
\ee
by the elements
\be
q\ot q'+q'\ot q, \qquad y\ot y'-y'\ot y, \qquad
q\ot y-y\ot q
\ee
for all $q,q'\in Q_x$, $y,y'\in Y_x$, and $x\in
X$. 
The $C^\infty(X)$-modules $\cA_S^k$ of sections of the
vector bundles
$\ol S^k\to X$ make up a direct system with respect to the
natural monomorphisms $\cA_S^k\to \cA_S^{k+1}$. Its
direct limit $\cA^\infty_S$ is endowed with a structure of
a graded commutative $C^\infty(X)$-ring generated by odd and
even elements. Generalizing the above technique for a graded
manifold $(X,\cA_Q)$ to $(X,\cA^\infty_S)$, one obtains the BGDA
isomorphic to $\cS^*_\infty$ \cite{book00,mpla}.
\end{rem}

\begin{rem} \label{g234} \mar{g234}
In physical applications, one can think
of $\cS^*_\infty$ as being a graded algebra of even and odd
variables on a smooth manifold $X$. In particular, this is the
case of the above mentioned antifield BRST theory on $X=\Bbb
R^n$
\cite{barn95,barn,bran,bran01}. Recall that, in gauge theory on
a principal bundle
$P\to X$ with a structure Lie group $G$, principal
connections on
$P\to X$ are represented by sections of the quotient 
$C=J^1P/G\to X$ \cite{book,book00,epr}. The connection
bundle
$C\to X$ is affine. It is coordinated by
$(x^\la, a^r_\la)$ such that, given a section $A$ of $C\to X$, 
its components $A^r_\la=a^r_\la\circ A$ are coefficients of
the familiar local connection form (i.e., gauge potentials).
Let $J^\infty C$ be the infinite order jet
space of $C\to X$ coordinated by
$(x^\la,a^r_{\la,\La})$, $0\leq |\La|$, and let
$\cP^*_\infty(C)$ be the polynomial subalgebra of the GDA
$\cO^*_\infty(C)$ of exterior forms of finite jet order on
$J^\infty C$ whose coefficients are polynomials of
$a^r_{\la,\La}$. Infinitesimal generators of one-parameter
groups of vertical automorphisms  (gauge transformations) of a
principal bundle $P$ are
$G$-invariant vertical vector fields on $P\to X$.
They are associated
to sections of the vector bundle
$V_GP=VP/G\to X$ of right Lie algebras of the group $G$. Let us
consider the simple graded manifold
$(X,\gA_{V_GY})$ constructed from this vector bundle. Its local
basis is $(x^\la, C^r)$. Let $\cC^*_{J^rQ}$ be the  
BGDA of graded
exterior forms on the graded manifold $(X,\gA_{J^rV_GP})$, and
$\cC^*_\infty(V_GP)$ the direct limit of the direct system
(\ref{g205}) of these algebras. Then the graded product
\mar{g211}\beq
\cS^*_\infty(V_G,C)=\cC^*_\infty(V_GP)\w \cP^*_\infty(C)
\label{g211}
\eeq
describe gauge potentials, odd ghosts and their jets in the
BRST theory. 
A generic basis for BRST theory contains the following
three sectors: (i) even classical fields of vanishing ghost
number (e.g., the above mentioned gauge potentials), (ii) odd
ghosts, ghosts-for-ghosts and antifields, (iii) even
ghosts-for-ghosts and antifields \cite{barn,bran01,gom}.
From the physical viewpoint, it seems more natural
to describe odd and even elements of the non-classical
sectors (ii) -- (iii) of this basis in the framework of the
unified construction in Remark \ref{g207} \cite{book00,mpla}.
\end{rem}

Hereafter, let the
collective symbols
$s^A_\La$ and $\th^A_\La$ stand both for even and odd
generating elements
$c^a_\La$, $y^i_\La$, $\th^a_\La$, $\th^i_\La$ of the
$C^\infty(X)$-algebra
$\cS^*_\infty$. Similarly to $\cO^*_\infty$, the BGDA
$\cS^*_\infty$ is decomposed  into
$\cS^0_\infty$-modules $\cS^{k,r}_\infty$ of
$k$-contact and
$r$-horizontal graded forms, together with the corresponding
projections $h_k$ and $h^r$. Accordingly, the graded exterior
differential
$d$ on 
$\cS^*_\infty$ is split into the sum $d=d_H+d_V$ of 
the total and vertical differentials 
\be
d_H(\f)=dx^\la\w d_\la(\f), \qquad d_V(\f)=\th^A_\La\w\dr^\La_A
\f, \qquad \f\in \cS^*_\infty.
\ee
The projection
endomorphism $\vr$ of $\cS^*_\infty$ is given by the expression
\be
\vr=\op\sum_{k>0} \frac1k\ol\vr\circ h_k\circ h^n,
\qquad \ol\vr(\f)= \op\sum_{|\La|\geq 0}
(-1)^{\nm\La}\th^A\w [d_\La(\dr^\La_A\rfloor\f)], 
\qquad \f\in \cS^{>0,n}_\infty,
\ee
similar to (\ref{r12}). The graded variational operator
$\dl=\vr\circ d$ is introduced. Then the BGDA
$\cS^*_\infty$ is split into the graded variational
bicomplex, analogous to the bicomplex (\ref{7}). 

The key point
is that, in contrast with the variational bicomplex (\ref{7}),
the algebraic Poincar\'e lemma has been stated only for the
short variational complex
\mar{g111}\beq
0\ar \Bbb R\ar
\cS^0_\infty\ar^{d_H}\cS^{0,1}_\infty \cdots
\ar^{d_H} \cS^{0,n}_\infty\ar^\dl 0, 
\label{g111}  
\eeq
of the BGDA $\cS^*_\infty$, i.e., this complex on $X=\Bbb R^n$
is exact at all terms, except $\Bbb R$
\cite{barn,drag}. We also consider the complex
\mar{g112}\beq
 0\to \cS^{1,0}_\infty\ar^{d_H} \cS^{1,1}_\infty\cdots
\ar^{d_H}\cS^{1,n}_\infty\ar^\vr E_1\to 0, \qquad E_1=
\vr(\cS^{1,n}_\infty), \label{g112}
\eeq
of graded one-forms, analogous to the complex (\ref{g201}), and
the de Rham complex
\mar{g110}\beq
0\to\Bbb R\ar \cS^0_\infty\ar^d \cS^1_\infty\cdots
\ar^d\cS^k_\infty \ar\cdots. \label{g110}
\eeq
One can think of elements 
\be
L=\cL\om\in \cS^{0,n}_\infty, \qquad 
\dl (L)= \op\sum_{|\La|\geq 0}
 (-1)^{|\La|}\th^A\w d_\La (\dr^\La_A L)\in E_1
\ee
of the complexes (\ref{g111}) -- (\ref{g112}) as being a graded
Lagrangian and its Euler--Lagrange operator, respectively.

\begin{theo} \label{g96} \mar{g96} The cohomology of
the complexes (\ref{g111}) and (\ref{g110}) equals the de
Rham cohomology $H^*(X)$ of 
$X$. The complex (\ref{g112}) is exact.
\end{theo}

\begin{proof} Next Section is devoted to the proof.
\end{proof}

\begin{cor} \label{cmp26} \mar{cmp26}
Every $d_H$-closed graded form $\f\in\cS^{0,m<n}_\infty$
falls into the sum
\mar{g214}\beq
\f=\vf + d_H\xi, \qquad \xi\in \cS^{0,m-1}_\infty, \label{g214}
\eeq 
where $\vf$ is a closed $m$-form on $X$. Every
$\dl$-closed graded form  (a variationally
trivial graded Lagrangian) $L\in \cS^{0,n}_\infty$ is the sum
\mar{g215}\beq
\f=\vf + d_H\xi, \qquad \xi\in \cS^{0,n-1}_\infty, \label{g215}
\eeq
where $\vf$ is a non-exact $n$-form on $X$. 
\end{cor}

The global exactness of the complex (\ref{g112}) at the term
$\cS^{1,n}_\infty$ results in the following.

\begin{prop} \label{g103} \mar{g103}
Given a graded Lagrangian $L$, there
is the decomposition 
\mar{g99,'}\ben
&& dL=\dl L - d_H(\Xi),
\qquad \Xi\in \cS^{1,n-1}_\infty, \label{g99}\\
&& \Xi=\op\sum_{s=0}
\th^A_{\nu_s\ldots\nu_1}\w
F^{\la\nu_s\ldots\nu_1}_A\om_\la,\qquad 
F_A^{\nu_k\ldots\nu_1}=
\dr_A^{\nu_k\ldots\nu_1}\cL-d_\la F_A^{\la\nu_k\ldots\nu_1}
+h_A^{\nu_k\ldots\nu_1},  \label{g99'}
\een
where local graded functions $h$ obey the relations
$h^\nu_a=0$,
$h_a^{(\nu_k\nu_{k-1})\ldots\nu_1}=0$. 
\end{prop}

\begin{proof} The proof of the decomposition (\ref{g99}) repeats
that in Proposition \ref{g93}. The coordinate expression
(\ref{g99'}) results from a direct computation.
\end{proof}

Proposition \ref{g103} states the existence of a global
finite order Lepagean equivalent 
$\Xi_L=\Xi+L$ of any graded Lagrangian $L$. Locally, one can
always choose $\Xi$ (\ref{g99'}) where all functions $h$ vanish.

\section{Proof of Theorem 4.1}

The proof of Theorem \ref{g96} follows a scheme of the proof
of Theorem \ref{g90}, but all sheaves are set on $X$.

 We start from the exactness of the
complexes (\ref{g111}) -- (\ref{g110}).
The exactness of the short variational complex (\ref{g111})
(the algebraic Poincar\'e lemma) and the de Rham complex
(\ref{g110}) (the Poincar\'e lemma)  on
$X=\Bbb R^n$ at all terms, except
$\Bbb R$, has been stated \cite{barn,bran,drag}. Let us
extend the algebraic Poincar\'e lemma to the complex
(\ref{g112}).  

\begin{lem} \label{g220} \mar{g220}
The complex (\ref{g112}) on $X=\Bbb R^n$ is exact.
\end{lem}

\begin{proof}
The fact that a $d_H$-closed graded form
$\f\in \cS^{1,m<n}_\infty$ is $d_H$-exact results from  
the algebraic Poincar\'e
lemma for horizontal graded forms 
$\f\in \cS^{0,m<n}_\infty$.  
Indeed, let us formally associate to a graded
$m$-form
$\f=\sum\f_A^\La\w ds^A_\La$ the horizontal graded
$(m-1)$-form $\ol\f=\sum\f_A^\La\ol s^A_\La$ depending on
additional variables $\ol s^A_\La$ of the same Grassmann parity
as
$s^A_\La$, and let us introduce the modified total differential
\be
\ol d_H\ol\f=d_H\ol\f + dx^\la\w (\ol s^A_\la\ol\dr_A\ol\f+
\ol s^A_{\la\m}\ol\dr_A^\m\ol\f+\cdots),
\qquad
\ol\dr^\La_A=\dr/\dr s^A_\La.
\ee
It is
easily justified that
$\ol{d_H\f}=\ol d_H\ol\f$. If $d_H\f=0$, then $\ol d_H\ol\f=0$
and, consequently, $\ol \f= \ol d_H \ol\psi$ where
$\ol\psi=\sum\psi_A^\La\ol s^A_\La$ is linear in $\ol
s^A_\La$.   Then $\f=d_H\psi$ where $\psi=\sum\psi_A^\La\w
ds^A_\La$. It remains to show that, if 
\be
\vr(\f)=\op\sum_{|\La|\geq 0}(-1)^{|\La|}\th^A\w
[d_\La(\dr_A^\La\rfloor\f)]= \op\sum_{|\La|\geq 0}
(-1)^{|\La|}\th^A\w
[d_\La\f_A^\La]=0, \qquad  \f\in
\cS^{1,n}_\infty,
\ee
then $\f$ is $d_H$-exact. A direct
computation gives 
\be
\f=d_H\psi,\qquad  \psi=-\op\sum_{|\La|\geq
0}\op\sum_{\Si+\Xi=\La} (-1)^{|\Si|}\th^A_{\Xi}\w
d_\Si\f^{\La+\m}_A \om_\m.
\ee
\end{proof}

Let us associate to each open subset $U\subset X$ the $\Bbb
R$-module $\cS^*_U$ of elements of the
$C^\infty(X)$-algebra $\cS^*_\infty$ whose coefficients
are restricted to
$U$. These modules make up a presheaf on $X$. Let
$\gS^*_\infty$ be the sheaf constructed from this presheaf and
$\G(\gS^*_\infty)$ its structure module of sections. One can
show that $\gS^*_\infty$ inherits the  bicomplex operations,
and $\G(\gS^*_\infty)$ does so. 
For short, we say that $\G(\gS^*_\infty)$ consists of
polynomials in $s^a_\La$, $ds^a_\La$ of
locally bounded jet order $|\La|$. There is the monomorphism
$\cS^*_\infty\to\G(\gS^*_\infty)$. 
Let us consider the
complexes of sheaves of $C^\infty(X)$-modules 
\mar{g114,5,3}\ben
&& 0\ar \Bbb R\ar
\gS^0_\infty\ar^{d_H}\gS^{0,1}_\infty \cdots
\ar^{d_H} \gS^{0,n}_\infty\ar^\dl 0, 
\label{g114}\\
&& 0\to \gS^{1,0}_\infty\ar^{d_H} \gS^{1,1}_\infty\cdots
\ar^{d_H}\gS^{1,n}_\infty\ar^\vr \gE_1\to 0, \qquad
\gE_1=\vr(\gS^{1,n}_\infty), 
\label{g115}\\
&& 0\to\Bbb R\ar \gS^0_\infty\ar^d \gS^1_\infty\cdots
\ar^d\gS^k_\infty \ar\cdots \label{g113}
\een
on $X$ and the complexes of their structure modules
\mar{g117,8,6}\ben
&& 0\ar \Bbb R\ar
\G(\gS^0_\infty)\ar^{d_H}\G(\gS^{0,1}_\infty) \cdots
\ar^{d_H} \G(\gS^{0,n}_\infty)\ar^\dl 0, 
\label{g117}\\
&& 0\to \G(\gS^{1,0}_\infty)\ar^{d_H} \G(\gS^{1,1}_\infty)\cdots
\ar^{d_H}\G(\gS^{1,n}_\infty)\ar^\vr \G(\gE_1)\to 0,
\label{g118}\\
&& 0\to\Bbb R\ar \G(\gS^0_\infty)\ar^d \G(\gS^1_\infty)\cdots
\ar^d\G(\gS^k_\infty) \ar\cdots. \label{g116} 
\een
The terms $\gS^{*,*}_\infty$ of these complexes are sheaves of
$C^\infty(X)$-modules. Therefore, they are fine and,
consequently, acyclic. Turn to the sheaf $\gE_1$.

\begin{lem} \label{g221} \mar{g221}
The sheaf $\gE_1$ on $X$ is fine.
\end{lem}

\begin{proof}
We use the fact that the sheaf $\gE_1$ is a projection
$\vr(\gS^{1,n}_\infty)$ of the sheaf $\gS^{1,n}_\infty$ of
$C^\infty(X)$-modules. Let $\gU=\{U_i\}_{i\in I}$ be a locally
finite open covering of $X$ and $\{f_i\in C^\infty(X)\}$ the
associated partition of unity. For any open subset and any
section $\f$ of the sheaf $\gS^{1,n}_\infty$ over $U$, let us
put $h_i(\f)=f_i\f$. The endomorphisms $h_i$ of
$\gS^{1,n}_\infty$ yield the $\Bbb R$-module endomorphisms
\be
\ol h_i=\vr\circ h_i: \gE_1\ar^{\rm in}
\gS^{1,n}_\infty\ar^{h_i}\gS^{1,n}_\infty\ar^\vr \gE_1
\ee 
of the sheaf $\gE_1$. They possess the properties for $\gE_1$
to be a fine sheaf. Indeed, for each $i\in I$, supp$f_i\subset
U_i$ provides a closed set such that $\ol h_i$ is zero outside
it, while the sum $\op\sum_{i\in I}\ol h_i$ is the
identity morphism.
\end{proof}

Consequently, the sheaf $\gE_1$ is acyclic.
The above mentioned exactness of the complexes (\ref{g111}) --
(\ref{g110}) on $X=\Bbb R^n$ implies the exactness of the
complexes of sheaves (\ref{g114}) -- (\ref{g113}) at all 
terms, except
$\Bbb R$.  It follows that the complexes (\ref{g114}) and 
(\ref{g113})  are resolutions of the constant sheaf $\Bbb R$,
while the complex (\ref{g115}) is exact.  By virtue of the
abstract de Rham theorem, the cohomology of the complexes
(\ref{g117}) and (\ref{g116}) equals the de Rham cohomology
$H^*(X)$ of $X$, whereas the complex (\ref{g118}) is
globally exact. It remains to prove the following.

\begin{theo} \label{g239} \mar{g239}
Cohomology of the complexes (\ref{g111})
-- (\ref{g110}) equals that of the complexes (\ref{g117}) --
(\ref{g116}). 
\end{theo}

The rest of this Section is the proof of Theorem \ref{g239}. Let
the common symbols
$\G^*_\infty$ and
$D$ stand for all the modules and the coboundary operators,
respectively, in the complexes (\ref{g117}) -- (\ref{g116}).
With this notation, one can say that any
$D$-closed element
$\f\in\G^*_\infty$ takes the form 
\mar{g230}\beq
\f=\vf + D\xi,  \label{g230}
\eeq
where $\vf$
is a non-exact closed exterior form on $X$. For the proof, it
suffices to show that, if an element
$\f\in
\cS^*_\infty$ is
$D$-exact in the module $\G^*_\infty$, then it is
so in $\cS^*_\infty$.
By virtue of the above mentioned Poincar\'e lemmas, 
if $X$ is contractible and a $D$-exact element $\f$ 
is of finite jet order
$[\f]$ (i.e., $\f\in\cS^*_\infty$), there exists an 
element $\vf\in
\cS^*_\infty$ such that $\f=D\vf$. Moreover,
a glance at the corresponding homotopy operators shows that 
the jet order
$[\vf]$ of $\vf$ is bounded by an integer $N([\f])$,
depending only on $[\f]$. We agree to call this fact the finite
exactness of an operator
$D$. Given an arbitrary manifold
$X$, the finite exactness takes place on 
any domain $U\subset X$. Let us state the following.

\begin{lem} \label{g222} \mar{g222}
Given a family $\{U_\al\}$ of disjoint open subsets of $X$,
let us suppose that the finite exactness of an operator $D$
takes place on every subset $U_\al$. Then it holds on the
union
$\op\cup_\al U_\al$.
\end{lem}

\begin{proof} Let
$\f\in\cS^*_\infty$ be a $D$-exact graded form on
$X$.
The finite exactness on $\cup
U_\al$ holds since $\f=D\varphi_\al$ on every
$U_\al$ and all $[\varphi_\al]<N([\f])$. 
\end{proof}

\begin{lem} \label{g223} \mar{g223}
Suppose that the finite exactness of an operator $D$ takes
place on open subsets
$U$, $V$ of $X$ and their non-empty overlap $U\cap V$. Then it
is also true on $U\cup V$.
\end{lem}

\begin{proof} Let
$\f=D\varphi\in\cS^*_\infty$ be a $D$-exact graded form on
$X$. By assumption, it can be brought into the form
$D\varphi_U$ on $U$ and $D\varphi_V$ on
$V$, where
$\varphi_U$ and $\varphi_V$ are graded forms of bounded
jet order.  Due to the decomposition (\ref{g230}), one can
choose the forms $\f_U$, $\f_V$ such that $\vf-\vf_U$ on $U$ and
$\vf-\vf_V$ on $V$ are $D$-exact. 
Let us consider the difference
$\varphi_U-\varphi_V$ on 
$U\cap V$.
It is a $D$-exact graded form of bounded
jet order which, by assumption, can be written as 
$\varphi_U-\varphi_V=D\si$ where 
$\si$ is also of bounded jet order. 
Lemma
\ref{am20} below shows that $\si=\si_U +\si_V$ where
$\si_U$ and
$\si_V$ are graded forms of bounded jet order on
$U$ and
$V$, respectively. Then, putting
\be
\varphi'_U=\varphi_U-D\si_U, \qquad  
\varphi'_V=\varphi_V+ D\si_V,
\ee
we have the graded form $\f$, equal to $D\varphi'_U$ on $U$ and 
$D\varphi'_V$ on $V$, respectively. Since the difference
$\vf'_U-\vf'_V$ on $U\cap V$ vanishes, we obtain $\f=D\vf'$ on
$U\cup V$  where 
\be
\vf'\op=^{\rm def}\left\{
\begin{array}{ccc}
\vf'|_U &= & \vf'_U\\
\vf'|_V &= & \vf'_V
\end{array}\right.
\ee
is of bounded jet order. 
\end{proof}

\begin{lem} \label{am20} \mar{am20}
Let $U$ and $V$ be open subsets of $X$ and $\si$
 a graded form of bounded jet order on
$U\cap V$. Then $\si$
splits into the sum $\si_U+ \si_V$ of graded exterior
forms $\si_U$ on
$U$ and
$\si_V$ on $V$ of bounded jet order. 
\end{lem} 

\begin{proof}
By taking a smooth partition of unity on $U\cup V$ subordinate 
to its cover
$\{U,V\}$ and passing to the function with support in $V$, 
we get a
smooth real function
$f$ on
$U\cup V$ which is 0 on a neighborhood $U_{U-V}$ of $U-V$ and 1
on a neighborhood $U_{V-U}$ of
$V-U$ in $U\cup V$. The graded form
$f\si$ vanishes on $U_{U-V}\cap (U\cap V)$ and, therefore, can
be extended by 0 to $U$. Let us denote it $\si_U$.
Accordingly, the graded form
$(1-f)\si$ has an extension $\si_V$ by 0 to 
$V$. Then $\si=\si_U +\si_V$ is a desired
decomposition because $\si_U$ and $\si_V$
are of finite jet order which does not exceed that of $\si$.
\end{proof}

By virtue of Lemmas \ref{g222} and \ref{g223}, the finite
exactness of an operator $D$ on a manifold $X$ takes place
because one can choose the corresponding cover of $X$ (see
Lemma 9.5 in \cite{bred}, Chapter V).

\section{Generalized Lagrangian supersymmetries}

A graded derivation  $\up\in\gd
\cS^0_\infty$ of the $\Bbb R$-ring $\cS^0_\infty$ is said to
be a generalized supersymmetry if the Lie derivative
$\bL_\up\f$ preserves the ideal of contact graded forms of the
BGDA
$\cS^*_\infty$.  

\begin{prop} \label{g231} \mar{g231} With respect to the local
basis $(x^\la,s^A_\La, dx^\la,\th^A_\La)$ for the BGDA
$\cS^*_\infty$, any generalized supersymmetry takes the form
\mar{g105}\beq
\up=\up_H+\up_V=\up^\la d_\la + (\up^A\dr_A +\op\sum_{|\La|>0}
d_\La\up^A\dr_A^\La), \label{g105}
\eeq
where $\up^\la$, $\up^a$ are local graded functions.
\end{prop}

\begin{proof}
The key point is that any element of the $C^\infty(X)$-algebra
$\cS^*_\infty$ is a section of a finite-dimensional vector
bundle over $X$ and any graded form is a finite
composition of $df$, $f\in\cS^0_\infty$. Therefore, the proof
follows those of Propositions \ref{g62} --
\ref{g72}.
\end{proof}

The interior product $\up\rfloor\f$ 
and the Lie derivative $\bL_\up\f$,
$\f\in\cS^*_\infty$ obey  
the same formulae
\be
&& \up\rfloor \f=\up^\la\f_\la + (-1)^{[\f_A]}\up^A\f_A, \qquad
\f\in \cS^1_\infty,\\ 
&& \up\rfloor(\f\w\si)=(\up\rfloor \f)\w\si
+(-1)^{|\f|+[\f][\up]}\f\w(\up\rfloor\si), \qquad \f,\si\in
\cS^*_\infty \\
&& \bL_\up\f=\up\rfloor d\f+ d(\up\rfloor\f), \qquad
\bL_\up(\f\w\si)=\bL_\up(\f)\w\si
+(-1)^{[\up][\f]}\f\w\bL_\up(\si).
\ee
as those on a graded manifold.
In particular, it
is easily justified that any vertical generalized supersymmetry
$\up$ (\ref{g105}) satisfies the relations 
\mar{g232,3,'}\ben
&& \up\rfloor d_H\f=-d_H(\up\rfloor\f), \label{g232}\\
&& \bL_\up(d_H\f)=d_H(\bL_\up\f), \qquad \f\in\cS^*_\infty.
\label{g233}
\een

\begin{prop}  \label{g106} \mar{g106}
Given a graded Lagrangian $L\in\cS^{0,n}_\infty$, its 
Lie derivative $\bL_\up L$
along a generalized
supersymmetry
$\up$ (\ref{g105}) fulfills the first variational formula
\mar{g107}\beq
\bL_\up L= \up_V\rfloor\dl L +d_H(h_0(\up\rfloor \Xi_L)) 
+ d_V (\up_H\rfloor\om)\cL, \label{g107}
\eeq
where $\Xi_L=\Xi+L$ is a Lepagean equivalent of $L$ given by
the coordinate expression (\ref{g99'}).
\end{prop}

\begin{proof} The proof follows that 
 of Proposition \ref{g75} and results from the decomposition
(\ref{g99}) and the relation (\ref{g232}). 
\end{proof}

In particular, let $\up$ be a divergence
symmetry of a graded Lagrangian $L$, i.e., $\bL_\up
L=d_H\si$,
$\si\in
\cS^{0,n-1}_\infty$. 
Then the first variational formula (\ref{g107}) 
restricted to Ker$\,\dl L$ leads to the weak
conservation law
\mar{g108}\beq
0\ap d_H(h_0(\up\rfloor\Xi_L)-\si). \label{g108}
\eeq
Similarly to the case of generalized symmetries, one can
justify that a generalized supersymmetry is a divergence
symmetry of a Lagrangian $L$ iff its vertical part $\up_V$
(\ref{g105}) is so. In this case the conservation law
(\ref{g108}) takes the form
\be
0\ap d_H(\up_V\rfloor\Xi_L-\si'),
\ee
where $\cJ_\up=\up_V\rfloor\Xi_L$ is the graded Noether current
along a generalized supersymmetry $\up_V$.

It should be emphasized that a Lepagean equivalent $\Xi_L$ in
the conservation law (\ref{g108}) is not uniquely 
defined. One can always choose it of the local form (\ref{g99'})
where graded functions $h$ vanish. There is a global
Lepagean equivalent of this form if either $\di X=1$ or $L$
is of first order.   

The BRST transformation in gauge theory on a principal bundle
$P\to X$ with a structure group $G$ in Remark \ref{g234}
exemplifies a vertical generalized supersymmetry. With respect
to a local basis
$(x^\la,a^r_\la, C^r)$ for the BGDA $\cS^*_\infty(V_G,C)$
(\ref{g211}),    it is given by the expression
\mar{g130}\ben
&& \up= \up_\la^r\frac{\dr}{\dr a_\la^r} +
\up^r\frac{\dr}{C^r}
+\op\sum_{|\La|>0}\left(d_\La\up_\la^r\frac{\dr}{\dr
a_{\La,\la}^r} + d_\La\up^r\frac{\dr}{C^r_\La}\right),
\label{g130}\\ 
&& \up_\la^r=C_\la^r +c^r_{pq}a^p_\la C^q,
\qquad \up^r=
\frac12c^r_{pq}C^p C^q, \nonumber 
\een
where $c^r_{pq}$ are structure constants of the right Lie
algebra of $G$. A remarkable peculiarity of this generalized
supersymmetry is that the
Lie derivative $\bL_\up$ along $\up$ (\ref{g130}) is
nilpotent on the module $S^{0,*}_\infty$ of horizontal
graded forms.

One says that a vertical generalized
supersymmetry
$\up$ (\ref{g105}) is nilpotent if 
\mar{g133}\beq
\bL_\up(\bL_\up\f)= \op\sum_{|\Si|\geq 0,|\La|\geq 0 }
(\up^B_\Si\dr^\Si_B(\up^A_\La)\dr^\La_A + 
(-1)^{[B][\up^a]}\up^B_\Si\up^A_\La\dr^\Si_B \dr^\La_A)\f=0
\label{g133}
\eeq
for any horizontal graded form $\f\in S^{0,*}_\infty$. A glance
at the second term in the expression (\ref{g133}) shows that a
nilpotent generalized supersymmetry is necessarily odd.
The following is an important criterion of a nilpotent
generalized supersymmetry. 

\begin{lem} A generalized supersymmetry $\up$ is 
nilpotent iff the equality
\be
\bL_\up(\up^A)=\op\sum_{|\Si|\geq 0} \up^B_\Si\dr^\Si_B(\up^A)=0
\ee
holds for all $\up^A$.
\end{lem}

\begin{proof}
The proof results from a direct computation.
\end{proof}

\begin{rem} \label{g235} \mar{g235}
A useful example of a nilpotent generalized
supersymmetry is the supersymmetry
\mar{g134}\beq
\up=\up^A(x)\dr_A +\op\sum_{|\La|>0}\dr_\La \up^A\dr_A^\La,
\label{g134} 
\eeq
where all $\up^A$ are smooth real functions on $X$, but
all
$s^A$ are odd.
\end{rem}

\section{Cohomology of nilpotent generalized supersymmetries}

Let $\up$ be a nilpotent generalized supersymmetry.
Since the Lie derivative $\bL_\up$ obeys the relation
(\ref{g233}), let us assume that the $\Bbb
R$-module $\cS^{0,*}_\infty$ of graded horizontal forms is
split into a bicomplex $\{S^{k,m}\}$ with respect
to the nilpotent operator $\bs_\up$ (\ref{g240}) and the total
differential
$d_H$.  This
bicomplex 
\be
d_H: S^{k,m}\to S^{k,m+1}, \qquad \bs_\up: S^{k,m}\to
S^{k+1,m} 
\ee
is graded by the form degree $0\leq m\leq
n$ and an integer $k\in \Bbb Z$, though it may happen that
$S^{k,*}=0$ starting from some number. For the sake
of brevity, let us call
$k$ the charge number. 

For instance, the BRST bicomplex
$S^{0,*}_\infty(C,V_GP)$ is graded by the charge number $k$
which is the polynomial degree of its elements 
 in odd variables $C^r_\La$. In this case, $\bs_\up$
(\ref{g240}) is the BRST operator. Since the ghosts $C^r_\La$
are characterized by the ghost number 1, $k$ is the ghost
number. The bicomplex defined by the supersymmetry (\ref{g134})
in Remark
\ref{g235} has the similar gradation, but the nilpotent
operator $\bs_\up$ decreases the odd polynomial
degree. 

Let us consider horizontal graded forms $\f\in \cS^{0,*}_\infty$
such that a nilpotent generalized supersymmetry $\up$ is their
divergence symmetry, i.e., $\bs_\up\f=d_H\si$. As was
mentioned above, we come to the relative and iterated cohomology
of the nilpotent operator
$\bs_\up$ (\ref{g240}) with respect to the total differential
$d_H$.

Recall that a horizontal graded form $\f\in S^{*,*}$ is said to
be  a relative (i.e., $(\bs_\up/d_H)$-) closed form if
$\bs_\up\f$ is a
$d_H$-exact form. This form is called
exact if it is a sum of an $\bs_\up$-exact form and a
$d_H$-exact form. Accordingly, we have the
relative cohomology $H^{*,*}(\bs_\up/d_H)$. If a
$(\bs_\up/d_H)$-closed form $\f$ is also $d_H$-closed, it is
called an iterated 
$(\bs_\up|d_H)$-closed
form. This form $\f$ is said to be exact if $\f=\bs_\up\xi
+d_H\si$, where $\xi$ is a $d_H$-closed form. Thus, we come
to the iterated cohomology 
$H^{*,*}(\bs_\up|d_H)$ of the
$(\bs_\up,d_H)$-bicomplex
$S^{*,*}$. It is the term $E_2^{*,*}$ of the
spectral sequence of this bicomplex
\cite{mcl}. There is an obvious isomorphism 
$H^{*,n,}(\bs_\up/d_H)=H^{*,n}(\bs_\up|d_H)$ of relative  and
iterated cohomology groups on horizontal graded densities. 
Forthcoming Theorems \ref{g250} and
\ref{g251} extend our results on iterated cohomology in
\cite{lmp} to an arbitrary nilpotent generalized supersymmetry. 

Note that, with respect to the total differential $d_H$, 
the bicomplex $S^{*,*}$ is the complex
\mar{g238}\beq
0\ar \Bbb R\ar
\cS^0_\infty\ar^{d_H}\cS^{0,1}_\infty \cdots
\ar^{d_H} \cS^{0,n}_\infty\ar^{d_H} 0,
\label{g238} 
\eeq 
which differs from
the short variational complex (\ref{g111}) in the last morphism.

\begin{prop} \label{g241} \mar{g241}
Cohomology groups $H^{m<n}(d_H)$ of the complex (\ref{g238})
equal  the de Rham cohomology groups $H^{m<n}(X)$ of $X$, while
the cohomology group $H^n(d_H)$ fulfills the
relation
\mar{g242}\beq
H^n(d_H)/H^n(X)=E_1. \label{g242}
\eeq
\end{prop}

\begin{proof} Cohomology of the complex (\ref{g238}) is 
determined similarly to that of the short variational complex 
(\ref{g111}). The
only difference is that the horizontal complex (\ref{g238}) on
$X=\Bbb R^n$ is not exact at the last term 
$\cS^{0,n}_\infty$. Accordingly, the corresponding complex of
sheaves 
\be
0\ar \Bbb R\ar
\gS^0_\infty\ar^{d_H}\gS^{0,1}_\infty \cdots
\ar^{d_H} \gS^{0,n}_\infty\ar^{d_H} 0
\ee
on $X$ fails to be a resolution of the constant sheaf $\Bbb R$
at the last term. Therefore, one should use a minor modification
of the abstract de Rham theorem \cite{jmp,tak2} in order
to obtain cohomology of the corresponding complex of structure
modules 
\be
0\ar \Bbb R\ar
\G(\gS^0_\infty)\ar^{d_H}\G(\gS^{0,1}_\infty) \cdots
\ar^{d_H} \G(\gS^{0,n}_\infty)\ar^{d_H} 0
\ee
at all the terms, except the last one. Then Theorem
\ref{g239} shows that this cohomology coincides with that of the
complex (\ref{g238}). The relation (\ref{g242}) results from
the formula (\ref{g215}). 
\end{proof}

\begin{rem}
Let us mention that, in  
the antifield BRST theory \cite{barn,bran} extended to an
arbitrary
$X$, the antibracket is defined on elements of the quotient
$H^n(d_H)/H^n(X)$. They correspond to local
functionals up to surface integrals. At the same time, the
antibracket on local functionals implies rather intricate
geometric interpretation of antifields \cite{khud,wit}, 
\end{rem}

\begin{theo} \label{g250} \mar{g250}
There is an epimorphism
\mar{g252}\beq
\zeta:
H^{m<n}(X)\to H^{*,m<n}(\bs_\up|d_H) \label{g252}
\eeq 
of the de Rham cohomology $H^m(X)$ of $X$ 
of form degree less than
$n$ onto the iterated
cohomology $H^{*,m<n}(\bs_\up|d_H)$.
\end{theo}

\begin{proof}
Since a nilpotent
generalized supersymmetry $\up$ is vertical, all
exterior forms
$\f$ on $X$ are $\bs_\up$-closed. It follows that they are
$(\bs_\up|d_H)$-closed. Since any $d_H$-exact horizontal graded
form is also $(\bs_\up|d_H)$-exact, we have     
a morphism $\zeta$ (\ref{g252}).
By virtue of Corollary \ref{cmp26} (and, equivalently,
Proposition \ref{g241}), any $d_H$-closed horizontal graded
$(m<n)$-form $\f$ is split into the sum $\f=\vf
+d_H\xi$ (\ref{g214}) of a closed $m$-form $\vf$ on $X$ and a
$d_H$-exact graded form. It follows that the morphism $\zeta$
(\ref{g252}) is an epimorphism.  
\end{proof}

In particular, if $X=\Bbb R^n$, the iterated cohomology
$H^{*,0<m<n}(\bs_\up|d_H)$ is trivial in contrast with the
relative ones.

The kernel of the morphism $\zeta$ (\ref{g252}) consists of
elements whose representatives are $\bs_\up$-exact
closed exterior forms on $X$. For instance, a glance
at the the BRST
transformation $\up$ (\ref{g130}) shows that, in BRST theory,
exterior forms on $X$ are never $\bs_\up$-exact. However, this
is not the case of the supersymmetry (\ref{g134}). 
Note that, in the both examples, exterior forms on $X$ are only
of zero charge number. In this case, we have the trivial 
iterated cohomology 
$H^{\neq 0,m<n}(\bs_\up|d_H)$ and an epimorphism (in
particular, an isomorphism) of $H^{m<n}(X)$ to
$H^{0,m<n}(\bs_\up|d_H)$.

Turn now to the iterated cohomology $H^{*,n}(\bs_\up|d_H)$.
It requires a particular analysis because, by virtue of
Proposition \ref{g241}, the cohomology $H^n(d_H)$ of the
complex (\ref{g238}) fails to equal the de Rham cohomology
$H^n(X)$ of $X$.

The bicomplex $S^{*,*}$ is
a complex with respect to the total coboundary operator
$\wt\bs_\up=\bs_\up+ d_H$.
 We aim
to state the relation between the iterated cohomology
$H^{*,m}(\bs_\up|d_H)$ and the total
$\wt\bs_\up$-cohomology $H^*(\wt\bs_\up)$ of the bicomplex
$S^{*,*}$. Similarly to the morphism (\ref{g252}),
there exists the morphism
\mar{g255}\beq
\g: H^{<n}(X)\to H^*(\wt\bs_\up) \label{g255}
\eeq
of the de Rham cohomology $H^{<n}(X)$ of $X$ 
of form degree 
$<n$ to the total
cohomology $H^*(\wt\bs_\up)$. Its kernel consists of
elements whose representatives are
$\wt\bs_\up$-exact closed exterior forms on $X$. Put $\ol H^*=
H^*(\wt\bs_\up)/\im \g$.

\begin{theo} \label{g251} \mar{g251}
 There is the
isomorphism 
\mar{g270}\beq
H^{*,n}(\bs_\up|d_H)/\ol H^*=\Ker \g. \label{g270}
\eeq
\end{theo}

\begin{proof} The proof falls into the following three steps.

(i) At first, we state the morphism
\mar{g271}\beq
\eta:H^{*,n}(\bs_\up|d_H)\to \Ker \g \label{g271}
\eeq
of the iterated cohomology $H^{*,n}(\bs_\up|d_H)$ to $\Ker\g$.
Let a horizontal graded $n$-form $\f_n$ be
$(\bs_\up|d_H)$-closed. Then, by definition,  
$\bs_\up\f_n$ is
$d_H$-exact, i.e.,
\mar{aa10}\beq
\bs_\up\f_n + d_H\f_{n-1}=0. \label{aa10}
\eeq
Acting on this equality by $\bs_\up$, we
observe that
$\bs_\up\f_{n-1}$ is a 
$d_H$-closed graded form, i.e.,
\mar{g140}\beq
\bs_\up\f_{n-1} + d_H\f_{n-2}=\vf_{n-1}, \label{g140}
\eeq
where $\vf_{n-1}$ is a closed $(n-1)$-form on $X$ in accordance
with Corollary
\ref{cmp26}. Since
$\bs_\up\vf_{n-1}=0$, an action of $\bs_\up$ on the equation
(\ref{g140}) shows that 
$\bs_\up\f_{n-2}$ is a 
$d_H$-closed graded form, i.e.,
\be
\bs_\up\f_{n-2} + d_H\f_{n-3}=\vf_{n-2}, 
\ee
where $\vf_{n-2}$ is a closed $(n-2)$-form on $X$.
Iterating the
arguments, one comes to the system of equations
\mar{aa11}\beq
\bs_\up\f_{n-k} + d_H\f_{n-k-1}=\vf_{n-k},
\qquad 0\leq k <n, \qquad
\bs_\up \f_0=\vf_0={\rm const},\label{aa11}
\eeq
which assemble into the descent equation 
\mar{g256,aa13}\ben
&& \wt\bs_\up\wt \f_n=\wt\vf_{n-1}, \label{g256} \\
&& \wt \f_n=\f_n+\f_{n-1}+\cdots +\f_0, \qquad \wt\vf_{n-1}=
\vf_{n-1}+\cdots +\vf_0.
\label{aa13}
\een
Thus, any $(\bs_\up|d_H)$-closed  horizontal graded form defines
a descent equation (\ref{g256}) whose right-hand side
$\wt\vf_{n-1}$ is  a closed exterior form on $X$ such that its
de Rham class belongs to the kernel
$\Ker\g$ of the morphism (\ref{g255}). 
For the sake
of brevity, let us denote this descent equation by
$\lng\wt\vf_{n-1}\rng$. Accordingly, we say that a
horizontal  graded form 
$\wt
\f_n$ (\ref{aa13}) is a solution of the descent equation
$\lng\wt\vf_{n-1}\rng$ (\ref{aa13}). A descent equation defined
by a $(\bs_\up|d_H)$-closed horizontal graded form $\f_n$ is
not unique. Let $\wt\f'$ be another solution of another descent
equation $\lng\wt\vf'_{n-1}\rng$
such that $\f_n=\f'_n$. Let us denote
$\Delta\f_k=\f_k-\f'_k$ and 
$\Delta\vf_k=\vf_k-\vf'_k$. Then the equation (\ref{aa10})
leads  to the equation $d_H(\Delta\f_{n-1})=0$. It follows that
\mar{g257}\beq
\Delta\f_{n-1}=d_H\xi_{n-2} +\al_{n-1}, \label{g257}
\eeq
where $\al_{n-1}$ is a closed $(n-1)$-form on $X$. Accordingly,
the equation (\ref{aa11}) leads to the equation
\be
\bs_\up(\Delta\f_{n-1})+d_H(\Delta\f_{n-2})=\Delta\vf_{n-1}.
\ee
Substituting the equality (\ref{g257}) into this equation, we
obtain the equality
\be
d_H(-\bs_\up\xi_{n-2}+ \Delta\f_{n-2})= \Delta\vf_{n-1}.
\ee
It follows that 
\be
\Delta\f_{n-2}
=\bs_\up\xi_{n-2} + d_H\xi_{n-3} +\al_{n-2}, \qquad
\Delta\vf_{n-1}=d\al_{n-2}
\ee
where $\al_{n-2}$ is an exterior form on $X$. Iterating the
arguments, one comes to the relations
\mar{g258}\beq
\Delta\f_{n-k} =\bs_\up\xi_{n-k} + d_H\xi_{n-k-1} + \al_{n-k},
\qquad
\Delta\vf_{n-k}=d \al_{n-k-1}, \qquad 1<k<n,
\label{g258}
\eeq
where $\al_{n-k-1}$ are exterior forms on $X$ and, finally, to
the equalities $\Delta\f_0=0$, $\Delta\vf_0=0$.
Then it
is easily justified that 
\mar{g259,'}\ben
&& \wt\f_n-\wt\f'_n=\wt\bs\wt\si +\wt\al, 
\qquad \wt\si= \Delta\f_{n-1}+\cdots +\Delta\f_1, \label{g259}\\
&& \wt\vf_{n-1} - \wt\vf'_{n-1}=d\wt\al, \qquad
\wt\al= \al_{n-1}+\cdots +\al_1. \label{g259'}
\een
It follows that right-hand sides of any two descent equations
defined by a $(\bs_\up|d_H)$-closed horizontal graded form
$\f_n$ differ from each other in an exact form on $X$. 
Moreover,  let $\f_n$ and $\f'_n$ be representatives of 
the same iterated cohomology class in $H^{*,n}(\bs_\up|d_H)$,
i.e., $\f_n=\f'_n+ \bs_\up\xi_n +d_H\si_{n-1}$. Let $\f_n$
provide a solution $\wt\f_n$ of a descent equation 
$\lng\wt\vf_{n-1}\rng$. Then $\f'_n$ defines a solution 
$\wt\f'=\wt\f+\wt\bs_\up(\xi_n+\si_{n-1})$ of the same descent
equation.
Thus, the assignment 
$\f_n\mapsto \lng\wt \vf_{n-1}\rng$ yields the desired
morphism $\eta$ (\ref{g271}).

(ii) Let $\wt\vf_{n-1}$ be a closed exterior form  on $X$ whose
de Rham cohomology class belongs to
$\Ker
\g$. Then $\wt\vf_{n-1}$ yields some descent equation
$\lng\wt\vf_{n-1}\rng$ (\ref{g256}). Let $\wt\vf'_{n-1}$ differ
from 
$\wt\vf_{n-1}$ in an exact form, i.e., let the relation
(\ref{g259'}) hold. Then any solution
$\wt
\f_n$ of the equation $\lng\wt\vf_{n-1}\rng$ yields a
solution
$\wt\f'_n=\wt\f_n-\wt\al$ (\ref{g259}) of the equation
$\lng\wt\vf'_{n-1}\rng$ such that $\f'_n=\f_n$.
It follows that the morphism $\eta$ (\ref{g271}) is an
epimorphism.

(iii) The kernel of the morphism $\eta$ (\ref{g271})
is represented by  $(\bs_\up|d_H)$-closed  horizontal
graded forms
$\f_n$ which define the homogeneous descent equation 
\mar{g273}\beq
\wt\bs_\up\wt\f_n=0. \label{g273}
\eeq
Its solutions $\wt\f_n$ are $\wt\bs_\up$-closed horizontal
graded forms. Let us assign to a solution $\wt\f_n$ of the
descent equation (\ref{g273}) its higher term $\f_n$. Running
back the arguments at the end of item (i), one can show that, if
solutions
$\wt\f_n$ and
$\wt\f'_n$ of the descent equation (\ref{g273}) belong to
the same total cohomology class,
then its higher terms $\f_n$ and
$\f'_n$ belong to the same iterated cohomology class. It
follows that the assignment $\wt\f_n\mapsto\f_n$ provides an
epimorphism of the total cohomology $H^*(\wt\bs_\up)$ onto 
$\Ker\eta$. The kernel of this epimorphism is represented by
solutions $\wt\f_n$ of the descent equation (\ref{g273})
whose higher term vanishes. Following item (i), one can
easily show that these solutions take the form
$\wt\f_n=\wt\bs_\up\si+\wt\al$, where $\wt\al$ is a closed
exterior form on $X$ of form degree $<n$. Cohomology
classes of these solutions exhaust the image of the
morphism $\g$ (\ref{g255}), i.e., $\im\g=\Ker\eta$.
\end{proof}

In particular, if the morphism $\g$ (\ref{g255}) is a
monomorphism (i.e., no non-exact closed exterior form on $X$ is 
$\wt\bs_\up$-exact), the isomorphism (\ref{g270}) gives the
isomorphism 
\be
H^{*,n}(\bs_\up|d_H)=H^*(\wt\bs_\up)/H^{<n}(X).
\ee
For instance, this is the case of the BRST transformation
(\ref{g130}) \cite{lmp}.

\end{document}